\documentclass{amsart}

\usepackage{amsmath}
\usepackage{amscd}
\usepackage{amsfonts}
\usepackage{amssymb}
\usepackage{latexsym}
\usepackage{pifont}
\usepackage{eepic}
\usepackage{array}

\newcommand{\ncm}{\newcommand}

\ncm{\beq}{\begin{equation}}
\ncm{\eeq}{\end{equation}}
\ncm{\bea}{\begin{eqnarray}}
\ncm{\eea}{\end{eqnarray}}
\ncm{\beanon}{\begin{eqnarray*}}
\ncm{\eeanon}{\end{eqnarray*}}

\newtheorem{thm}{Theorem}[section]
\newtheorem{pro}[thm]{Proposition}
\newtheorem{lem}[thm]{Lemma}
\newtheorem{cor}[thm]{Corollary}

\theoremstyle{definition}
\newtheorem{defi}[thm]{Definition}

\theoremstyle{remark}
\newtheorem{rmk}[thm]{Remark}

\numberwithin{equation}{section}

\def\Set{\mathsf{Set}}

\def\M{\mathsf{M}}

\def\A{\mathcal{A}}
\def\B{\mathcal{B}}
\def\C{\mathcal{C}}
\def\D{\mathcal{D}}

\ncm{\F}{\mathcal{F}}
\ncm{\G}{\mathcal{G}}

\ncm{\V}{\mathcal{V}}
\ncm{\W}{\mathcal{W}}

\ncm{\asso}{\mathbf{a}}
\ncm{\luni}{\mathbf{l}}
\ncm{\runi}{\mathbf{r}}

\ncm{\End}{\operatorname{End}}
\def\Hom{\mbox{\rm Hom}\,}

\def\id{\mbox{\rm id}\,}
\def\Center{\mbox{\rm Center}\,}
\ncm{\Ind}{\operatorname{Ind}}

\newcommand{\ci}{\circ}
\newcommand{\bo}{\,\Box\,}
\def\o{\otimes}
\def\x{\times}

\ncm{\amalgo}[1]{\underset{\scriptscriptstyle #1}{\o}}
\ncm{\oT}{\amalgo{T}}
\ncm{\oB}{\amalgo{B}}
\ncm{\oR}{\amalgo{R}}
\ncm{\oL}{\amalgo{L}}
\ncm{\ex}[1]{\underset{\scriptscriptstyle #1}{\x}}
\def\bra{\langle}
\def\ket{\rangle}
\def\under{\mbox{\rm\_}\,}
\ncm{\rarr}[1]{\stackrel{#1}{\longrightarrow}}
\ncm{\larr}[1]{\stackrel{#1}{\longleftarrow}}
\def\cop{\Delta}

\def\eps{\varepsilon}
\ncm{\op}{\mathrm{op}}
\ncm{\coop}{\mathrm{coop}}
\ncm{\Fi}{\varphi}
\ncm{\IN}{\mathrm{in}}
\ncm{\OUT}{\mathrm{out}}
\ncm{\OR}{\overrightarrow}
\ncm{\OL}{\overleftarrow}

\def\oneT{^{(1)}}
\def\twoT{^{(2)}}
\def\threeT{^{(3)}}
\def\oneB{_{(1)}}
\def\twoB{_{(2)}}
\def\threeB{_{(3)}}
\def\fourB{_{(4)}}
\def\oneR{^{[1]}}
\def\twoR{^{[2]}}

\def\oneL{_{[1]}}
\def\twoL{_{[2]}}

\def\PL{\pi_{\scriptscriptstyle L}}
\def\PR{\pi_{\scriptscriptstyle R}}
\ncm{\I}{\mathcal{I}}
\def\la{\!\rightharpoonup\!}
\def\ra{\!\leftharpoonup\!}
\def\du1{\hat 1}
\def\iso{\rarr{\sim}}

\def\under{\mbox{\rm\_}\,}
\ncm{\Cnt}{\mathsf{C}}
\ncm{\ZZ}{\mathbb{Z}}

\ncm{\GRB}{\Gamma_{\scriptscriptstyle RB}}
\ncm{\GBR}{\Gamma_{\scriptscriptstyle BR}}
\ncm{\GLT}{\Gamma_{\scriptscriptstyle LT}}
\ncm{\GTL}{\Gamma_{\scriptscriptstyle TL}}
\ncm{\GLB}{\Gamma_{\scriptscriptstyle LB}}
\ncm{\GBL}{\Gamma_{\scriptscriptstyle BL}}
\ncm{\GRT}{\Gamma_{\scriptscriptstyle RT}}
\ncm{\GTR}{\Gamma_{\scriptscriptstyle TR}}

\ncm{\LA}{\mathcal{L}}
\ncm{\BA}{\mathcal{B}}
\ncm{\RA}{\mathcal{R}}
\ncm{\TA}{\mathcal{T}}

\ncm{\Oo}{\mathcal{O}}
\ncm{\Ha}{\mathcal{H}}
\ncm{\Ve}{\mathcal{V}}

\ncm{\Mnull}{\mathbf{M}(\G^0,\varnothing)}

\ncm{\rhomb}[1]{
\begin{picture}(20,20)
\put(10,0){\line(-1,1){10}}
\put(10,0){\line(1,1){10}}
\put(10,20){\line(-1,-1){10}}
\put(10,20){\line(1,-1){10}}
\put(8,7){$#1$}
\end{picture}
}

\begin{document}

\title{The double algebraic view of finite quantum groupoids}
\author[K. Szlach\'anyi]{Korn\'el Szlach\'anyi}
\address{Research Institute for Particle and Nuclear Physics, Budapest\\
H-1525 Budapest, P. O. Box 49, Hungary}
\email{szlach@rmki.kfki.hu}

\begin{abstract}
Double algebra is the structure modelled by the properties of the ordinary and the convolution product in Hopf algebras, weak Hopf algebras and Hopf algebroids if a Frobenius integral is given. The Hopf algebroids possessing a Frobenius integral are precisely the Frobenius double algebras in which the two multiplications satisfy distributivity. The double algebra approach makes it manifest that all comultiplications in such measured Hopf algebroids are of the Abrams-Kadison type, i.e., they come from a Frobenius algebra structure in some bimodule category.
Antipodes for double algebras correspond to the Connes-Moscovici `deformed' antipode as we show by discussing Hopf and weak Hopf algebras from the double algebraic point of view. Frobenius algebra extensions provide further examples that need not be distributive.
\end{abstract}

\thanks{Partially supported by the Hungarian Scientific Research Fund, OTKA T 034 512 and T 043 159.}

\maketitle

\vskip 2truecm

\section{Introduction}

Let $A$  be a Hopf algebra, weak Hopf algebra or a Hopf algebroid \cite{HGD}
and assume that there exists a (left or right) integral $i\in A$ which is a Frobenius homomorphism on the dual algebra. Such `measured quantum groupoids' are known to have two algebra structures:
the underlying algebra, which we call the vertical algebra $V$, and the
horizontal algebra $H$ with multiplication given by the convolution product $\star$ and with unit given by the integral $i$. Studying the interrelation of these two algebras leads to the following

\begin{defi} \label{def: DA}
Let $k$ be a commutative ring. A $k$-module $A$ equipped with two associative unital $k$-algebra structures $V=\bra A,\ci,e\ket$ and $H=\bra A,\star,i\ket$ is called a double algebra over $k$ if the following properties hold:
\begin{quote}
\begin{itemize}
\item[\bf A1.] $(a\star e)\ci b = ((a\star e)\ci i)\star b$
\item[\bf A2.] $a\ci(b\star e) = (i\ci(b\star e))\star a$
\item[\bf A3.] $(a\ci i)\star b = ((a\ci i)\star e)\ci b$
\item[\bf A4.] $a\star(b\ci i) = (e\star(b\ci i))\ci a$
\item[\bf A5.] $a\ci(e\star b) = a\star(i\ci(e\star b))$
\item[\bf A6.] $(e\star a)\ci b = b\star((e\star a)\ci i)$
\item[\bf A7.] $a\star(i\ci b) = a\ci(e\star(i\ci b))$
\item[\bf A8.] $(i\ci a)\star b = b\ci((i\ci a)\star e)$
\end{itemize}
\end{quote}
for all $a,b\in A$.
\end{defi}
It will turn out to be useful to view the axioms in terms of the $k$-linear maps $A\to A$ defined by
\begin{alignat}{2}
\Fi_L(a)&:=a\star e &\qquad\Fi_R(a)&:=e\star a \label{eq: Fi1}\\
\Fi_B(a)&:=a\ci i &\qquad \Fi_T(a)&:=i\ci a \label{eq: Fi2}
\end{alignat}
In case of the quantum groupoids mentioned above these $\Fi$'s
are Frobenius homomorphisms onto subalgebras $L$ and $R$ of $V$ and $B$ and $T$ of $H$, respectively. $L$ and $R$ are the target and source subalgebras, traditionally called $A^L$ and $A^R$ \cite{BNSz} or $A_t$ and $A_s$ \cite{NV} of a weak Hopf algebra. If $i$ was chosen to be a right integral then $T$ is the trivial right $A$-module with cyclic vector $i$ and $B$ is the space of right integrals.

It is known from M. M\"uger's work \cite{Muger} that the left regular module in the monoidal category of left $H$-modules over a Frobenius Hopf algebra $H$ is a Frobenius algebra in which the multiplication is given by the convolution product.
The ordinary multiplication of $H$, however, does not belong to this category. The double algebra is just the structure that incorporates both multiplications in a
completely symmetric way. Moreover, it goes beyond Hopf algebras as far as the base
ring(s) need not be commutative.

The two comultiplications of the quantum groupoid arise naturally from the dual bases of $\Fi_B$ and $\Fi_T$ just as one associates a comultiplication to a Frobenius extension \cite{Kadison} or to a Frobenius algebra \cite{Abrams}.
This gives a convenient formalism to deal with Hopf algebroids because the many bimodule structures over $L$ and $R$ one needs in \cite{HGD} can all be replaced with one of the natural bimodules $_BA_B$, $_TA_T$, $_LA_L$, $_RA_R$ that arise from the subalgebras $B,T\subset H$ and $L, R\subset V$. Since the comultiplications are uniquely determined by the Frobenius homomorphisms, the measured quantum groupoids consist of nothing more than two Frobenius algebra structures with certain compatibility conditions.

The Frobenius double algebras are similar to the \textit{double Frobenius algebras} of M. Koppinen \cite{Koppinen} and to the \textit{bi-Frobenius algebras} of Y. Doi and M. Takeuchi \cite{Doy-Takeuchi} in which, however, the base algebras $B$, $L$, $T$, $R$ all coincide with the ground field.
The closest relatives of double algebras in the C$^*$-algebraic framework are probably the \textit{(quantum) hypergroups} \cite{Chap-Vainer}.

The use of four base algebras instead of two reveals a $D_4$ dihedral symmetry in measured quantum groupoids, which is evident from the above axioms.
Therefore, as a rough picture of the double algebra, we may represent $A$ as a square with its four boundary edges on the left, right, bottom and top corresponding to the base ideals.
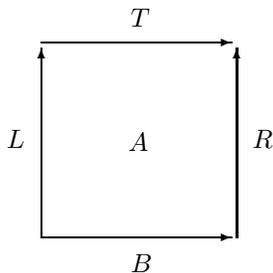
\[
\begin{picture}(100,100)
\put(13,13){\vector(1,0){72}} \put(47,0){$B$}
\put(13,13){\vector(0,1){72}} \put(0,47){$L$}
\put(87,13){\vector(0,1){72}} \put(93,47){$R$}
\put(13,87){\vector(1,0){72}} \put(47,93){$T$}
\put(46,46){$A$}
\end{picture}
\]
The orientation of the edges correspond to our convention of writing the second factor $b$ in a vertical multiplication $a\ci b$ on the top of $a$ and in a horizontal multiplication $a\star b$ on the right of $a$. Involving two neighbour base homomorphisms $\Fi$, each of the eight axioms can be associated to a corner.
The picture suggests a relation to double categories which is probably not accidental.
The examples of weak Hopf algebras constructed from double groupoids by N. Andruskiewitsch and S. Natale \cite{AN} also point to that direction.

Even if the four $\Fi$'s are Frobenius homomorphisms the double algebra is far from being a Hopf algebroid. The missing property can be most easily captured by saying that the two multiplications should be \textit{distributive} over each other (see Definition \ref{def: D}).
In Theorem \ref{thm: main} we prove that the Hopf algebroids with Frobenius integrals are precisely the \textit{distributive Frobenius double algebras}. Since any statements about a double algebra remains true if horizontal and vertical are interchanged, the dual Hopf algebroid appears to be built into the double algebra as well as the original: they are the horizontal and vertical Hopf algebroids of the double algebra. The arising picture is reminiscent to that of the 'double triangle algebras' \cite{Ocneanu,PetZub,CoqTri}.

Unlike the original papers \cite{Bohm,HGD,FrobD2} that are based on bialgebroids, the present double algebra approach to (measured) Hopf algebroids has a strong flavor of weak Hopf algebras, although the base algebras are no longer separable. This can be best seen in Section \ref{s: 2} or in the Maschke theorem of Section \ref{s: G}. The story from double algebras to Hopf algebroids is almost entirely contained in Sections \ref{s: 2}, \ref{s: F2} and \ref{s: D}. Section \ref{s: G} deals with an intermediate situation when distributivity does not hold but a Maschke theorem already works. Section \ref{s: S} serves for an introduction to the (double algebraic) antipode.

As for the possible significance of nondistributive double algebras the example of Subsection \ref{ss: F ext} is worth a mention. The two-step centralizer $\Cnt_{M_2}(N)$ in the Jones tower $N\subset M\subset M_2\subset M_3$ over any Frobenius extension $N\subset M$ is a double algebra with antipode. The convolution product of this double algebra is, of course, obtained from the algebra structure of the second two-step centralizer $\Cnt_{M_3}(M)$ by using Fourier transformation. In general this double algebra is not a quantum groupoid unless the extension $N\subset M$ is of depth 2.

\section{Double algebras} \label{s: 2}

If $A$ and $A'$ are double algebras in the sense of Definition \ref{def: DA} then a map of double algebras $f\colon A\to A'$ is simply a $k$-module map which is an algebra homomorphism both vertically and horizontally. Relaxed morphisms, such as partly nonunital ones, can also play a role but we will not need them here.

For any double algebra we define the maps $\Fi_B$, $\Fi_L$, $\Fi_R$ and $\Fi_T$ by equations (\ref{eq: Fi1}) and (\ref{eq: Fi2}) and rewrite the axioms in terms of them as follows.
\begin{quote}
\begin{itemize}
\item[\bf A1.] $\Fi_L(a)\ci b=\Fi_B\Fi_L(a)\star b$
\item[\bf A2.] $a\ci\Fi_L(b)=\Fi_T\Fi_L(b)\star a$
\item[\bf A3.] $\Fi_B(a)\star b=\Fi_L\Fi_B(a)\ci b$
\item[\bf A4.] $a\star\Fi_B(b)=\Fi_R\Fi_B(b)\ci a$
\item[\bf A5.] $a\ci\Fi_R(b)=a\star\Fi_T\Fi_R(b)$
\item[\bf A6.] $\Fi_R(a)\ci b=b\star\Fi_B\Fi_R(a)$
\item[\bf A7.] $a\star\Fi_T(b)=a\ci\Fi_R\Fi_T(b)$
\item[\bf A8.] $\Fi_T(a)\star b=b\ci\Fi_L\Fi_T(a)$
\end{itemize}
\end{quote}
(Note that we employed juxtaposition to denote composition of maps because $\ci$ is reserved for vertical multiplication.)

These expressions contain neither $i$'s nor $e$'s, so they can be the starting point of
further generalizations. But even in the unital case they are simpler to deal with than the original Definition \ref{def: DA}, at least after acquainting with the basic properties of the $\Fi$'s.

The $D_4$ symmetry of this structure can be generated by the following operations:
1. interchanging the vertical and horizontal algebra structures and 2. replacing $\ci$ with its opposite while keeping $\star$ unchanged.
Applying the first to the double algebra $A$ one obtains the double algebra $A^D$, applying the second one obtains $A_\op$. The double algebra $((A^D)_\op)^D$ is called $A_\coop$.

\begin{lem}
In any double algebra $A$ over $k$
\begin{itemize}
\item[\bf (a)] $L:=\Fi_L(A)$ and $R:=\Fi_R(A)$ are left, resp. right, ideals in $H$,
\item[\bf (b)] $B:=\Fi_B(A)$ and $T:=\Fi_T(A)$ are left, resp. right, ideals in $V$,
\item[\bf (c)] $L$ and $R$ are subalgebras of $V$,
\item[\bf (d)] $B$ and $T$ are subalgebras of $H$,
\item[\bf (e)] with respect to the natural bimodule structure for a subalgebra the $\Fi$'s are bimodule maps
\begin{alignat*}{2}
\Fi_L&\colon \,_LV_L\to \,_LL_L,&\qquad \Fi_R&\colon\,_RV_R\to\,_RR_R,\\
\Fi_B&\colon \,_BH_B\to \,_BB_B,&\qquad \Fi_T&\colon\,_TH_T\to\,_TT_T\,.
\end{alignat*}
\end{itemize}
\end{lem}
\begin{proof}
\textbf{(a)} and \textbf{(b)} are obvious.
Using axiom {\bf A1} twice and associativity of $\star$
\begin{equation*}
\Fi_L(\Fi_L(a)\ci b)=((a\star e)\ci i)\star(b\star e) =\Fi_L(a)\ci\Fi_L(b)\,,
\end{equation*}
together with $\Fi_L(i)=i\star e=e$, proves that $L$ is a (unital) subalgebra of $V$ and at the same time that $\Fi_L$ is a left $L$-module map. That it is also a right $L$-module map follows by applying {\bf A2}, then associativity of $\star$ and
then {\bf A2} again,
\begin{align*}
\Fi_L(a\ci\Fi_L(b))&=\Fi_L(\Fi_T\Fi_L(b)\star a)=\Fi_T\Fi_L(b)\star\Fi_L(a)=\\
&=\Fi_L(a)\ci\Fi_L(b)\,.
\end{align*}
Passing to $A_\op$ the above result implies that $R$ is also a subalgebra of $V$ and that $\Fi_R$ is a bimodule map. Passing to $A^D$ we obtain the corresponding results for $\Fi_B$ and $\Fi_T$.
\end{proof}

The $L$, $R$, $B$, and $T$ will be called respectively the left, right, bottom and top subalgebra, or ideal, of $A$. Together they will be referred to as the base ideals, or as the base algebras, and the $\Fi$'s as the base homomorphisms of the double algebra $A$.

\begin{lem} \label{lem: 2}
In any double algebra
\begin{itemize}
\item[\bf (a)] the base homomorphisms satisfy the identities
\begin{align*}
\Fi_X\Fi_Y\Fi_X\ =\ \Fi_X \quad \text{for }(X,Y)=&(L,B),(B,L),(B,R),(R,B),\\
&(R,T),(T,R),(T,L),\text{and }(L,T)\,,
\end{align*}
i.e., for any pair of base ideals which share the same corner, and
\[
\Fi_X\Fi_Y=\Fi_Y\Fi_X \qquad \text{for } (X,Y)=(L,R)\text{ or }(B,T)\,,
\]
\item[\bf (b)] restrictions of base maps give rise to algebra isomorphisms
\begin{alignat*}{2}
\Fi_L|_B\colon&B\iso L\quad&\text{with inverse}\quad\Fi_B|_L&:L\iso B\,,\\
\Fi_B|_R\colon&R\iso B^\op\quad&\text{with inverse}\quad\Fi_R|_B&:B^\op\iso R\,,\\
\Fi_R|_T\colon&T\iso R\quad&\text{with inverse}\quad\Fi_T|_R&:R\iso T\,,\\
\Fi_T|_L\colon&L\iso T^\op\quad&\text{with inverse}\quad\Fi_L|_T&:T^\op\iso L\,,
\end{alignat*}
\item[\bf (c)] and the $L$ and $R$ commute within $V$ and the $B$ and $T$ commute within $H$, i.e.,
\begin{align*}
l\ci r&=r\ci l\qquad l\in L,\ r\in R\,,\\
b\star t&=t\star b \qquad b\in B,\ t\in T\,.
\end{align*}
\end{itemize}
\end{lem}
\begin{proof}
(a) Setting $b=e$ in {\bf A1} one immediately gets $\Fi_L\Fi_B\Fi_L=\Fi_L$ and similarly, each further axiom provides one more identity.
That the maps under (b) are $k$-module isomorphisms is obvious from (a).
It remains to show that they are algebra maps. We are content with proving this for $\Fi_L|_B$.
\[
\Fi_L(b)\ci\Fi_L(b')=\Fi_B\Fi_L(b)\star\Fi_L(b')=b\star\Fi_L(b')=\Fi_L(b\star b').
\]
In order to prove (c) we compute
\begin{align*}
\Fi_L(a)\ci\Fi_R(b)&=((a\star e)\ci i)\star(e\star b)=a\star e\star b\\
\Fi_R(b)\ci\Fi_L(a)&=(a\star e)\star((e\star b)\ci i)=a\star e\star b
\end{align*}
where {\bf A1} was used in the first and {\bf A6} in the second line.
The dual formula yields commutativity of $T$ and $B$.
\end{proof}

\begin{cor} \label{cor: centers}
The appropriate restrictions of the base homomorphisms provide the algebra isomorphisms
\begin{align}
L\cap\Center V&\cong B\cap T\cong R\cap\Center V\\
B\cap\Center H&\cong L\cap R\cong T\cap\Center H\\
L\cap R\cap\Center V&\cong B\cap T\cap\Center H\,.
\end{align}
\end{cor}
In case of Hopf algebras (Subsection \ref{ss: H}) the base algebras are all trivial, i.e., coincide with $k\cdot e$ and $k\cdot i$, respectively. In case of weak Hopf algebras (Subsection \ref{ss: WHA}) the base algebras are separable $k$-algebras but all the commutative algebras in the Corollary can be nontrivial. For Hopf algebroids
(Subsection \ref{ss: HGD}) the base algebras are unrestricted. As in \cite{NV} one may call the situation $B\cap T=k\cdot i$ a \textit{connected double algebra}, the situation $L\cap R=k\cdot e$ a \textit{coconnected double algebra} and if both conditions are met a \textit{biconnected double algebra}.

The isomorphisms between the base subalgebras offers the following interpretation of the axioms. The eight possible actions of the four base algebras on the double algebra
are organized into four actions. Therefore we redraw the picture of the double algebra as
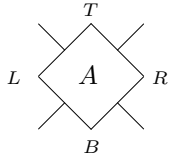
\[
\parbox[c]{2in}{
\begin{picture}(100,60)
\put(50,15){\line(-1,1){20}}
\put(50,15){\line(1,1){20}}
\put(50,55){\line(-1,-1){20}}
\put(50,55){\line(1,-1){20}}
\put(45,32){$A$}
\put(47,6){$\scriptstyle B$}
\put(18,32){$\scriptstyle L$}
\put(73,32){$\scriptstyle R$}
\put(47,58){$\scriptstyle T$}
\put(40,25){\line(-1,-1){10}}
\put(60,25){\line(1,-1){10}}
\put(40,45){\line(-1,1){10}}
\put(60,45){\line(1,1){10}}
\end{picture}
}
\]
where the attached lines correspond to the above mentioned four base algebra actions, but not to specific base algebras, however. The Temperley-Lieb type of identities under Lemma \ref{lem: 2} (a) suggest that the tensor squares $A\amalgo{X} A$ w.r.t.
$B$, $L$, $T$, $R$ should be represented respectively by the diagrams
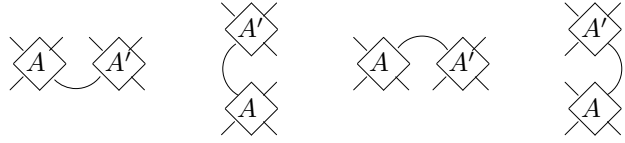
\[
\parbox[c]{4in}{
\begin{picture}(260,60)
\put(0,17){
\begin{picture}(70,40)
\put(20,5){\rhomb{\!\!A}}
\put(50,5){\rhomb{\!\!A'}}
\put(47,17){\arc{23}{0.65}{2.35}}
\put(68,10){\line(1,-1){5}}
\put(57,20){\line(-1,1){5}}
\put(27,20){\line(-1,1){5}}
\put(67,20){\line(1,1){5}}
\put(37,20){\line(1,1){5}}
\put(27,10){\line(-1,-1){5}}
\end{picture}
}
\put(70,0){
\begin{picture}(60,60)
\put(30,5){\rhomb{\!\!A}}
\put(30,35){\rhomb{\!\!A'}}
\put(44,30){\arc{23}{2.1}{4.1}}
\put(48,20){\line(1,1){5}}
\put(37,50){\line(-1,1){5}}
\put(37,10){\line(-1,-1){5}}
\put(48,40){\line(1,-1){5}}
\put(48,10){\line(1,-1){5}}
\put(48,50){\line(1,1){5}}
\end{picture}
}
\put(130,17){
\begin{picture}(70,40)
\put(20,5){\rhomb{\!\!A}}
\put(50,5){\rhomb{\!\!A'}}
\put(48,14){\arc{23}{3.8}{5.7}}
\put(68,10){\line(1,-1){5}}
\put(57,10){\line(-1,-1){5}}
\put(27,20){\line(-1,1){5}}
\put(67,20){\line(1,1){5}}
\put(38,10){\line(1,-1){5}}
\put(27,10){\line(-1,-1){5}}
\end{picture}
}
\put(200,0){
\begin{picture}(60,60)
\put(30,5){\rhomb{\!\!A}}
\put(30,35){\rhomb{\!\!A'}}
\put(42,30){\arc{23}{5.3}{7.3}}
\put(37,20){\line(-1,1){5}}
\put(37,50){\line(-1,1){5}}
\put(37,10){\line(-1,-1){5}}
\put(37,40){\line(-1,-1){5}}
\put(48,10){\line(1,-1){5}}
\put(48,50){\line(1,1){5}}
\end{picture}
}
\end{picture}
}
\]
where $A$ and $A'$ refer to the first and second tensorands, respectively.

The (a) part of the above Lemma implies that $B$, $L$, $T$, $R$ are direct summands of the $k$-module $A$. The projections $\Fi_L\Fi_B$, $\Fi_R\Fi_B$, \dots etc onto the base ideals appear also in the formulas
\begin{align*}
a\ci i&=\Fi_B(a)=\Fi_B(a)\star i=\Fi_L\Fi_B(a)\ci i\\
i\ci a&=\Fi_T(a)=i\star\Fi_T(a)=i\ci\Fi_R\Fi_T(a)
\end{align*}
which suggest that $i$ should be a 2-sided integral in $V$ in some appropriate quantum groupoid sense. The dual formulas present $e$ as a 2-sided integral in $H$. Although double algebras are far from being quantum groupoids, we shall use the name integral for the elements of the $k$-modules defined in the next Lemma. Since a double algebra unifies  two dual structures, it should not be a surprise that the integrals do not give entirely new ideals in $A$, just give a new characterization of the base ideals. More precisely, we have

\begin{lem} \label{lem: integrals}
In any double algebra $A$ define the $k$-submodules
\begin{align*}
\I_R&:=\{l\in A\,|\,l\star a=l\star\Fi_B\Fi_R(a),\ a\in A\,\} \\
\I_L&:=\{r\in A\,|\,a\star r=\Fi_B\Fi_L(a)\star r,\ a\in A\,\} \\
\I_T&:=\{b\in A\,|\,b\ci a=b\ci\Fi_L\Fi_T(a),\ a\in A\,\} \\
\I_B&:=\{t\in A\,|\,a\ci t=\Fi_L\Fi_B(a)\ci t,\ a\in A\,\}
\end{align*}
Then
\begin{alignat*}{2}
&L\subset\I_R\subset\Cnt_V(R)&\qquad&R\subset\I_L\subset\Cnt_V(L)\\
&B\subset\I_T\subset\Cnt_H(T)&\qquad&T\subset\I_B\subset\Cnt_H(B)
\end{alignat*}
If the bilinear forms $\Fi_{B/T}(\under\star\under)$ and $\Fi_{L/R}(\under\ci\under)$
are nondegenerate  then
\begin{equation*}
\I_R=L,\quad\I_L=R,\quad\I_T=B,\quad\I_B=T.
\end{equation*}
\end{lem}
\begin{proof}
It suffices to prove the statements for $\I_R$. In view of Lemma \ref{lem: 2} (a) for all $a'\in A$
\[
\Fi_L(a')\star a=a'\star \Fi_R(a)=a'\star\Fi_R\Fi_B\Fi_R(a)
=\Fi_L(a')\star\Fi_B\Fi_R(a)
\]
therefore $L\subset \I_R$. Next let $l'\in\I_R$. Then for all $a\in A$
\[
l'\ci\Fi_R(a)=l'\star\Fi_T\Fi_R(a)=l'\star\Fi_B\Fi_R(a)=\Fi_R(a)\ci l'
\]
hence $\I_R\subset\Cnt_V(R)$. Now assume that $\Fi_B(\under\star\under)$ is nondegenerate and let $l'\in \I_R$. Then
\begin{align*}
\Fi_B(l'\star a)&=\Fi_B(l'\star\Fi_B\Fi_R(a))=\Fi_B(l')\star\Fi_B(e\star a)=
\Fi_B(\Fi_B(l')\star e\star a)\\
&=\Fi_B(\Fi_L\Fi_B(l')\star a)
\end{align*}
implies that $l'\in L$.
\end{proof}

\begin{lem} \label{lem: inv}
An element $r\in R$ is invertible in $V$ iff it is invertible in $R$.
Similar statements hold for $L$, $B$ and $T$.
\end{lem}
\begin{proof}
Let $v\in A$ be such that $v\ci r=e=r\ci v$. Then
\begin{align*}
\Fi_R\Fi_T(v)\ci r&=\Fi_R(\Fi_T(v)\ci r)=\Fi_R\Fi_T(v\ci r)=\Fi_R\Fi_T(e)=e\\
r\ci\Fi_R\Fi_B(v)&=\Fi_R(r\ci\Fi_B(v)=\Fi_R\Fi_B(r\ci v)=\Fi_R\Fi_B(e)=e
\end{align*}
therefore $\Fi_R\Fi_T(v)=\Fi_R\Fi_B(v)=r^{-1}\in R$.
\end{proof}

Our next theme is the restriction to the base ideals of the would-be Nakayama automorphism of the base maps.
Composing the algebra isomorphisms found in Lemma \ref{lem: 2} we obtain two isomorphisms from $L$ to $R^\op$ and two ones from $B$ to $T^\op$
\begin{align*}
\Fi_R\Fi_T|_L:L\iso R^\op&\qquad\Fi_R\Fi_B|_L:L\iso R^\op\\
\Fi_T\Fi_L|_B:B\iso T^\op&\qquad\Fi_T\Fi_R|_B:B\iso T^\op
\end{align*}
the differences of which being measured by the automorphisms
\begin{equation*}
\Fi_R\Fi_B\Fi_L\Fi_T|_R\colon R\to R\,,\qquad
\Fi_T\Fi_R\Fi_B\Fi_L|_T\colon T\to T\,.
\end{equation*}

\begin{lem} \label{lem: Nakayama}
For any $a\in A$, $l\in L$, $b\in B$, $r\in R$ and $t\in T$ we have
\begin{alignat*}2
\Fi_L(a\ci r)&=\Fi_L(\nu_L(r)\ci a) &\quad\Fi_R(a\ci l)&=\Fi_R(\nu_R(l)\ci a)\\
\Fi_B(a\star t)&=\Fi_B(\nu_B(t)\star a) &\quad\Fi_T(a\star b)&=
\Fi_T(\nu_T(b)\star a)
\end{alignat*}
where we introduced
\begin{alignat}{2}
\nu_L(r)&=\Fi_R\Fi_B\Fi_L\Fi_T(r)\,,&\quad \nu_R(l)&=\Fi_L\Fi_B\Fi_R\Fi_T(l)\\
\nu_B(t)&=\Fi_T\Fi_L\Fi_B\Fi_R(t)\,,&\quad\nu_T(b)&=\Fi_B\Fi_L\Fi_T\Fi_R(b)\,.
\end{alignat}
\end{lem}
\begin{proof}
We prove the formula for $\Fi_B$, the rest is left to the reader.
\begin{align*}
\Fi_B(\Fi_T\Fi_L\Fi_B\Fi_R(t)\star a)&=\Fi_B(a\ci\Fi_L\Fi_B\Fi_R(t))=\\
&=a\ci\Fi_B(\Fi_L\Fi_B\Fi_R(t))=a\ci\Fi_B\Fi_R(t)=\Fi_B(a\ci\Fi_R(t))=\\
&=\Fi_B(a\star t)\,.
\end{align*}
\end{proof}

\section{Frobenius double algebras} \label{s: F2}

We recall some facts about Frobenius extensions \cite{Kadison} in order to fix the terminology.
For a subalgebra $L\subset A$ a Frobenius homomorphism $\Fi$ is a bimodule map $_LA_L\to\,_LL_L$ possessing dual bases. The latter means a finite set of pairs $\{x_i,y_i\}$ of elements of $A$ such that
\[
\sum_i \Fi(ax_i)y_i\ =\ a\ =\ \sum_i x_i\Fi(y_i a)\,,\quad \forall a\in A\,.
\]
The element $\sum_i x_i \o y_i\in A\oL A$ is independent of the choice of the dual bases because for any choice it is the unit for the associative multiplication defined on $A\oL A$ by means of $\Fi$ as $(a\o b)(c\o d)=a\Fi(bc)\o d$. For this reason we shall, by an abuse of language, call $\sum_i x_i\o y_i$ \textit{the} dual basis of $\Fi$. The Nakayama automorphism of a Frobenius homomorphism $\Fi$ is an automorphism of the centralizer $\Cnt_A(L)$ defined by either one of the equivalent equations
\begin{align}
\label{eq: Nakayama}
\nu_\Fi(c)&=\sum_i \Fi(x_i c) y_i\\
\label{eq: Nakayama psi}
\Fi(ac)=\Fi(\nu_\Fi(c) a)\,,&\quad a\in A,\ c\in\Cnt_A(L)\\
\label{eq: Nakayama u-v}
\sum_i x_ic \o y_i\ =\ \sum_i x_i\o\nu_\Fi(c)y_i\,,&\quad c\in\Cnt_A(L)\,.
\end{align}
The central element $\Ind\Fi:=\sum_ix_i y_i\in\Center A$ defines the index of $\Fi$
but not of the extension, however.

\begin{defi}
The double algebra $\bra A,\ci,\star,e,i\ket$ is called Frobenius if
$\Fi_L$, $\Fi_B$, $\Fi_R$ and $\Fi_T$ are Frobenius homomorphisms.
\end{defi}
We introduce a special notation for the dual bases of each Frobenius homomorphism,
\begin{itemize}
\item[] $\sum_k u_k\o v_k\in A\oB A$ is the dual basis of $\Fi_B$,
\item[] $\sum_j x_j\o y_j\in A\oL A$ is the dual basis of $\Fi_L$,
\item[] $\sum_k u^k\o v^k\in A\oT A$ is the dual basis of $\Fi_T$,
\item[] $\sum_j x^j\o y^j\in A\oR A$ is the dual basis of $\Fi_R$,
\end{itemize}
although the summation symbol will always be suppressed.

As we have the inclusion $T\subset\Cnt_A(B)$, the Nakayama automorphism $\nu_B$ of $\Fi_B$ can be restricted to $T$. Lemma \ref{lem: Nakayama} yields explicit expressions for this restriction and for the analogous ones of $\nu_L$, $\nu_T$ and $\nu_R$.

In order to see the connection of double algebras to quantum groupoids it is crucial
to change the view of Frobenius structures as just Frobenius homomorphisms.
As it has been made clear by L. Kadison \cite{Kadison} for Frobenius extensions and by L. Abrams \cite{Abrams} for Frobenius algebras a Frobenius homomorphism $A\to X$ for a subalgebra $X\subset A$ is equivalent to a comonoid in the bimodule category $_X\M_X$, i.e., an $X$-coring that is compatible with multiplication in the sense of the comultiplication being an $A$-$A$-bimodule map. Therefore in a Frobenius double algebra
\begin{itemize}
\item[] $\bra A, \Delta_B,\Fi_B\ket$ is a comonoid in $_B\M_B$, where
      $\Delta_B(a)=a\star u_k\oB v_k$,
\item[] $\bra A, \Delta_L,\Fi_L\ket$ is a comonoid in $_L\M_L$, where
      $\Delta_L(a)=a\ci x_j\oL y_j$,
\item[] $\bra A, \Delta_T,\Fi_T\ket$ is a comonoid in $_T\M_T$, where
      $\Delta_T(a)=a\star u^k\oT v^k$,
\item[] $\bra A, \Delta_R,\Fi_R\ket$ is a comonoid in $_R\M_R$, where
      $\Delta_R(a)=a\ci x^j\oR y^j$.
\end{itemize}
Compatibility of $\cop_B$, for example, with multiplication means that $a\star u_k\oB v_k=u_k\oB v_k\star a$, a well known property of the dual basis.
Although the compatibility conditions of $\cop_B$ with $\star$ are very different from what one needs in a bialgebra or in any quantum groupoid, this is not so with $\Delta_B$ and $\ci$.

\begin{pro} \label{pro: almost QGD}
Let $A$ be a Frobenius double algebra. Then both\footnote{The 6-tuple notation
compresses the total algebra, the base algebra the source map, the target map, the comultiplication and the counit, in this order.}
\[
\bra V,B,\Fi_L|_B,\Fi_R|_B,\Delta_B,\Fi_B\ket \quad\text{and}\quad
\bra H,L,\Fi_B|_L,\Fi_T|_L,\Delta_L,\Fi_L\ket
\]
satisfy the axioms of a left bialgebroid \cite{KSz}, except multiplicativity of the comultiplication and both
\[
\bra V,T,\Fi_R|_T,\Fi_L|_T,\Delta_T,\Fi_T\ket \quad\text{and}\quad
\bra H,R,\Fi_T|_R,\Fi_B|_R,\Delta_R,\Fi_R\ket
\]
satisfy the axioms of a right bialgebroid \cite{KSz}, except multiplicativity of the comultiplication.
\end{pro}
\begin{proof}
We prove the statement for the bialgebroid $V$ over $B$. The source and target maps $s_B:=\Fi_L|_B\colon B\to V$ and $t_B:=\Fi_R|_B\colon B^\op\to V$, respectively, are algebra maps the ranges of which commute by (b) and (c) part of Lemma \ref{lem: 2}. The corresponding $B$-$B$-bimodule structure on $V$,
$b\cdot a\cdot b'=s_B(b)\ci t_B(b')\ci a$, coincides with the natural bimodule structure $_BH_B$ via the equality $H=V=A$ as $k$-modules. Therefore the above comonoid structure $\bra A, \Delta_B,\Fi_B\ket$ is precisely the one that is needed for $V$ to be a left bialgebroid over $B$. The comultiplication and the counit preserve the unit because
\begin{gather*}
\Delta_B(e)=e\star u_k\o v_k=\Fi_R(u_k)\o v_k=e\star\Fi_B\Fi_R(u_k)\o v_k=\\
=e\o \Fi_B(e\star u_k)\star v_k=e\o e \\
\Fi_B(e)=e\ci i =i
\end{gather*}
The counit axioms (i.e. axiom (vii) on p. 80 of \cite{KSz}) now take the form
\[
\Fi_B(a\ci\Fi_L\Fi_B(a'))=\Fi_B(a\ci a')=\Fi_B(a\ci\Fi_R\Fi_B(a'))\,\quad
a,a'\in A
\]
and hold true because of the identities $\Fi_B(a\ci a'')=a\ci\Fi_B(a'')$ and $\Fi_B\Fi_L\Fi_B=\Fi_B=\Fi_B\Fi_R\Fi_B$. It remains to show the Takeuchi property
of the comultiplication, namely
\begin{equation} \label{eq: Takeuchi prop}
\Delta_B(a)\ci (\Fi_R(b)\o 1)=\Delta_B(a)\ci(1\o\Fi_L(b))\,\qquad a\in A,\ b\in B\,.
\end{equation}
Insert here the expression of $\Delta_B$ through the dual basis $u_k\o v_k$ and use
{\bf A2} and {\bf A5} to rewrite the statement as
\[
a\star u_k\star\Fi_T\Fi_R(b)\oB v_k=a\star u_k\oB \Fi_T\Fi_L(b)\star v_k\,.
\]
Now equation (\ref{eq: Nakayama u-v}) and the expression of $\nu_B|_T$ in Lemma \ref{lem: Nakayama} reduces the statement to proving that
\[
\Fi_T\Fi_L(b)=\Fi_T\Fi_L\Fi_B\Fi_R\Fi_T\Fi_R(b).
\]
But this plainly follows by repeatedly applying Lemma \ref{lem: 2} (a).
\end{proof}

Notice that the above proof explains the appearence of the Takeuchi property merely from the Frobenius structure of the double algebra. In other words, the Abrams-Kadison comultiplication automatically satisfies the Takeuchi property within a double
algebra. But it is not necessarily multiplicative as the example in Subsection \ref{ss: F ext} shows.

\section{The Galois maps} \label{s: G}

For $A$ a Frobenius double algebra we can define the maps, for the time being in $\M_k$,
\begin{align}
\Gamma_{XY}:\ A\amalgo{X}A\to A\amalgo{Y}A&\text{ where }X,Y\in\{L,B,R,T\}\text{ are neighbours} \notag\\
\GRB(a\oR a')&=a\star u_k\oB v_k\ci a' \label{eq: GRB}\\
\GBR(a\oB a')&=a\star x^j\oR y^j\ci a'\\
\GLT(a\oL a')&=a\ci u^k\oT v^k\star a'\\
\GTL(a\oT a')&=a\ci x_j\oL y_j\star a'\\
\GLB(a\oL a')&=u_k\ci a'\oB v_k\star a \label{eq: GLB}\\
\GBL(a\oB a')&=x_j\star a'\oL y_j\ci a\\
\GRT(a\oR a')&=a'\star u^k\oT a\ci v^k\\
\GTR(a\oT a')&=a'\ci x^j\oR a\star y^j
\end{align}
In order to see that they are well defined it suffices to show this for $\GRB$.
Using \textbf{A5}, centrality of the dual basis and \textbf{A5} again we obtain for $r\in R$
\begin{gather*}
(a\ci r)\star u_k\oB v_k\ci a'=a\star\Fi_T(r)\star u_k\oB v_k\ci a'=\\
a\star u_k\oB(v_k\star\Fi_T(r))\ci a'=a\star u_k \oB v_k\ci r\ci a'\,.
\end{gather*}
Also notice that due to that $a\star\under\in \End A_B$ and $\under \ci a'\in \End\,_BA$ the definition of $\GRB$ is independent of the choice of $u_k$, $v_k$ within the dual basis.

These maps will be called the Galois maps of the Frobenius double algebra because they all are variations of the formula $(\id\o \mu)(\Delta\o\id)$ with some multiplication $\mu$ and some comultiplication $\Delta$. The question is what module structures the $\Gamma$'s preserve? The $\GRB$ is an $H$-$V$-bimodule map in the obvious sense,
\[
h\star\GRB(a\oR a')\ci v:=h\star a\star u_k\oB v_k\ci a'\ci v
=\GRB(h\star a\oR a'\ci v)
\]
for $h,v\in A$. So is the $\GBR$. However, there are more interesting module structures to preserve. For example, $A\oR A$ is also a left $V$ module via
$v\ci(a\oR a')=v\ci a\oR a'$. Unfortunately $A\oB A$ does not have a left $V$-action unless we use the comultiplication $\Delta_B$ to define it. Roughly speaking preservation of left $V$-action by $\GRB$ requires multiplicativity of $\Delta_B$.
This property will not hold until Section \ref{s: D} so here we are content with considering $\GRB$ and $\GBR$ as $H$-$V$-bimodule maps, $\GLT$ and $\GTL$ as $V$-$H$-bimodule maps, $\GLB$ and $\GBL$ as right $H\o_k V$-module maps and $\GRT$ and $\GTR$ as left $V\o_k H$-module maps.
\begin{lem} \label{lem: G}
Let $A$ be a Frobenius double algebra. Then the Galois maps $\Gamma_{XY}$ are ivertible with inverse $\Gamma_{YX}$ iff the equations
\begin{align}
u_k\star(v_k\ci a)&=\Fi_T\Fi_R(a)  \label{eq: G1}\\
(a\star x^j)\ci y^j&=\Fi_L\Fi_B(a) \label{eq: G2}
\end{align}
\begin{align}
(a\ci u^k)\star v^k&=\Fi_B\Fi_L(a)\\
x_j\ci(y_j\star a)&=\Fi_R\Fi_T(a)\\
(u_k\ci a)\star v_k&=\Fi_T\Fi_L(a) \label{eq: G5}\\
(x_j\star a)\ci y_j&=\Fi_R\Fi_B(a)\\
u^k\star(a\ci v^k)&=\Fi_B\Fi_R(a)\\
x^j\ci(a\star y^j)&=\Fi_L\Fi_T(a)
\end{align}
are satisfied for all $a\in A$.
\end{lem}
\begin{proof}
It suffices to prove that (\ref{eq: G1}) is equivalent to $\GBR\GRB=\id$. Since $i\o e$ is a cyclic vector of $A\oR A$ as an $H$-$V$-bimodule,
\[
\id_{A\oR A}=\GBR\GRB\quad\Longleftrightarrow\quad i\oR e=u_k\star x^j\oR y^j\ci v_k\,.
\]
Using centrality of the dual basis of $\Fi_R$ and nondegeneracy of $\Fi_R$ this is equivalent to the validity, for all $a\in A$, of the equation
\[
u_k\star(v_k\ci x^j)\ci\Fi_R(y^j\ci a)=i\ci\Fi_R(e\ci a)
\]
which is the same as (\ref{eq: G1}) because $u_k\star\under$ is a right $R$-module map.
\end{proof}

For a double algebra which satisfies the conditions of the above Lemma the index of
a base homomorphism is calculable as
\begin{alignat}{2}
\Ind\Fi_L&=x_j\ci y_j&=\Fi_R(i\ci i) \quad& \in R\cap\Center V\\
\Ind\Fi_R&=x^j\ci y^j&=\Fi_L(i\ci i) \quad& \in L\cap\Center V\\
\Ind\Fi_B&=u_k\star v_k&=\Fi_T(e\star e)\quad& \in T\cap \Center H\\
\Ind\Fi_T&=u^k\star v^k&=\Fi_B(e\star e)\quad& \in B\cap\Center H\,.
\end{alignat}
By Lemma \ref{lem: inv} and Corollary \ref{cor: centers} invertibility of $\Ind\Fi_L$, for example, is equivalent to $i\ci i$ being invertible in $B\cap T$.
Consider the special case when $A$ is also biconnected. Then
invertibility of both $\Ind\Fi_B$ and $\Ind\Fi_L$ means that $u_k\star v_k=i\cdot\beta$ and $\Fi_B(i)=i\cdot \beta'$ for some units $\beta,\beta'\in k_\x$.
That is to say $\bra A,\star,i,\cop_B,\Fi_B\ket$ is a special Frobenius algebra in $_B\M_B$ in the sense of \cite{Fuchs-Schweigert}.

The next result is a Maschke type theorem for double algebras.
\begin{thm}
Let $\bra A,\ci,e,\star,i\ket$ be a Frobenius double algebra for which the conditions of
Lemma \ref{lem: G} are satisfied. Let $\M_V$ and $_V\M$ denote the category of left, respectively right $V$-modules. Then the following are equivalent:
\begin{enumerate}
\item The inclusion $B\hookrightarrow A$ is split mono in $_V\M$.
\item $\Fi_B\colon A\to B$ is split epi in $_V\M$.
\item $L\subset V$ is a separable extension of $k$-algebras.
\item $i$ is von Neumann regular element in $V$, i.e., there exists $j\in A$ s.t.
\newline $i\ci j \ci i=i$.
\item $R\subset V$ is a separable extension of $k$-algebras.
\item $\Fi_T\colon A\to T$ is split epi in $\M_V$.
\item The inclusion $T\hookrightarrow A$ is split mono in $\M_V$.
\end{enumerate}
Similar equivalences hold for $H$ ($e$ is regular in $H$ iff $B\subset H$ is separable iff\dots etc).
\end{thm}
\begin{proof}
Since $B$ and $T$ are left, resp. right principal ideals of $V$ generated by $i$,
the equivalence of (1), (2), (4), (6) and (7) is a well-known result in ring theory (see e.g. \cite[p. 175]{AF}). Assume (4). Then
\[
\Fi_L\Fi_B(i\ci j)=\Fi_L(i\ci j\ci i)=\Fi_L(i)=e\,.
\]
One may notice that this is a formula showing that the bottom integral $i\ci j$ is normalized in the sense of weak Hopf algebras \cite{BNSz}. Next we show that the $j$ can be chosen in $L$. As a matter of fact, let $l:=\Fi_L\Fi_T(j)$. Then $i\ci l=\Fi_T\Fi_L\Fi_T(j)=\Fi_T(j)=i\ci j$. Therefore the Galois property (\ref{eq: G2})
implies that
\[
e=\Fi_L\Fi_B(i\ci l)=((i\ci l)\star x^j)\ci y^j=(\Fi_T(l)\star x^j)\ci y^j=
x^j\ci l\ci y^j
\]
This means that the map
\[
\sigma_R\colon V\to A\oR V,\qquad a\mapsto a\star(x^j\ci l)\oR y^j
\]
is a $V$-$V$ bimodule map splitting the epimorphism
\[
\mu_R\colon V\oR V \to A,\qquad a\oR a'\mapsto a\ci a'
\]
defined by multiplication. This proves (4) $\Rightarrow$ (5).
Now assume (5) and let $e_k\oR f_k$ be a separating idempotent. Then for $q:=e_k\ci\Fi_R(f_k)\in\Cnt_V(R)$ we have
\begin{align*}
i\ci q\ci i&=e_k\ci\Fi_R(f_k\ci i)\ci i=e_k\ci\Fi_B\Fi_R\Fi_B(f_k)=e_k\ci\Fi_B(f_k)
=e_k\ci f_k\ci i\\
&=i
\end{align*}
This proves (5) $\Rightarrow$ (4). The equivalence (3) $\Leftrightarrow$ (4) can be seen analogously: the $j$ can be chosen to be $r\in R$ and then $\sigma_L(a)=a\ci x_j\ci r\oL y_j$ defines a splitting map for the multiplication $\mu_L\colon V\oL V\to V$. Vice versa, if $e_k\oL f_k$ is a separating idempotent for $L\subset V$ then $i$ is regular with $j$ equal to $e_k\ci\Fi_L(f_k)$.
\end{proof}

The extension $L\subset V$ is called split if there exists an $L$-$L$-bimodule map $\eps\colon V\to L$ such that $\eps(e)=e$. Since $L\subset V$ is a Frobenius
extension, it is split exactly when there exists an $r\in\Cnt_V(L)$ such that $r\star e=e$. In this case $\eps(a)=\Fi_L(r\ci a)$. If $e$ is a regular element of $H$, so $e=e\star j\star e$, then $r$ can be chosen to be $e\star j\in R$.
Now assume that the double algebra have Galois maps as in Lemma \ref{lem: G}.
If both $e$ and $i$ are regular in $H$ and $V$, respectively, then all the algebra extensions $L\subset V$, $R\subset V$, $B\subset H$ and $T\subset H$ are split separable Frobenius extensions.

\section{The antipode} \label{s: S}

\begin{defi} \label{def: S}
An antipode for the double algebra $\bra A,\ci,e,\star,i\ket$ is a $k$-module map
$S\colon A\to A$ such that
\begin{align}
\Fi_B(a'\star(a''\ci a))&=\Fi_B((a'\ci S(a))\star a'') \label{eq: SB}\\
\Fi_R(a'\ci(a\star a''))&=\Fi_R((S(a)\star a')\ci a'') \label{eq: SR}\\
\Fi_L((a'\star a)\ci a'')&=\Fi_L(a'\ci(a''\star S(a))) \label{eq: SL}\\
\Fi_T((a\ci a')\star a'')&=\Fi_T(a'\star(S(a)\ci a'')).\label{eq: ST}
\end{align}
\end{defi}

If the base homomorphisms are nondegenerate and antipode exists then it is a unique anti algebra endomorphism of both $H$ and $V$. In the rest of the Section we restrict ourselves to study antipodes in Frobenius double algebras.

For any double algebra $A$ the $k$-module $A$ carries four (left or right $H$ or $V$) actions given by
\begin{alignat}{2}
\TA_a&:=\under \ci a\in\End\,_BA_B&\qquad \RA_a&:=\under\star a\in\End\,_LA_L\\
\LA_a&:=a\star\under\in\End\,_RA_R&\qquad \BA_a&:=a\ci\under\in\End\,_TA_T
\end{alignat}
for $a\in A$. If $A$ is Frobenius then each of these \textit{regular actions} can be left or right transposed w.r.t. the appropriate Frobenius homomorphisms. Therefore, we define $\TA_a^<,\TA_a^>\in\End\,_BA_B$,\dots etc by the following formulae.
\begin{alignat*}{2}
\Fi_B(\TA_a^<(a')\star a'')&=\Fi_B(a'\star \TA_a(a''))&\quad
\Fi_B(a'\star \TA_a^>(a''))&=\Fi_B(\TA_a(a')\star a'')\\
\Fi_L(\RA_a^<(a')\ci a'')&=\Fi_L(a'\ci \RA_a(a''))&\quad
\Fi_L(a'\ci \RA_a^>(a''))&=\Fi_L(\RA_a(a')\ci a'')\\
\Fi_R(\LA_a^<(a')\ci a'')&=\Fi_R(a'\ci \LA_a(a''))&\quad
\Fi_R(a'\ci \LA_a^>(a''))&=\Fi_R(\LA_a(a')\ci a'')\\
\Fi_T(\BA_a^<(a')\star a'')&=\Fi_T(a'\star \BA_a(a''))&\quad
\Fi_T(a'\star \BA_a^>(a''))&=\Fi_T(\BA_a(a')\star a'')
\end{alignat*}

\begin{lem} \label{lem: corner S}
For any Frobenius double algebra the bimodule maps $\TA^<_\bullet(e)\colon a\mapsto \TA^<_a(e)$, \dots etc are invertible with the following inverses:
\begin{alignat}{2}
\left[\TA_\bullet^<(e)\right]^{-1}&=\LA_\bullet^>(i)&\qquad
\left[\LA_\bullet^<(i)\right]^{-1}&=\BA_\bullet^<(e)\\
\left[\BA_\bullet^>(e)\right]^{-1}&=\RA_\bullet^<(i)&\qquad
\left[\RA_\bullet^>(i)\right]^{-1}&=\TA_\bullet^>(e)
\end{alignat}
\end{lem}
\begin{proof}
Any $a\in A$ can be expressed in the following two ways:
\begin{align*}
a&=u_k\star\Fi_B(v_k\star a)=\Fi_R\Fi_B(v_k\star a)\ci u_k=\Fi_R(\LA_{v_k}(a)\ci i)
\ci u_k\\
&=\Fi_R(a\ci \LA_{v_k}^>(i))\ci u_k\\
a&=\Fi_R(a\ci x^j)\ci y^j=y^j\star\Fi_B\Fi_R(a\ci x^j)=y^j\star\Fi_B(e\star \TA_{x^j}(a))\\
&=y^j\star\Fi_B(\TA^<_{x^j}(e)\star a)
\end{align*}
The first implies that $\LA_{v_k}^>(i)\o u_k$ is a right unit in the unital algebra $A\oR A$ therefore
\[
\LA_{v_k}^>(i)\oR u_k=x^j\oR y^j\,.
\]
Similarly, the second expression implies that
\[
y^j\oB \TA_{x^j}^<(e)=u_k\oB v_k\,.
\]
This means precisely that $\LA^>_\bullet(i)$ is the inverse of $\TA^<_\bullet(e)$. Employing the symmetries of the double algebra axioms the remaining three
relations are consequences.
\end{proof}
Notice that the eight maps in the Lemma are related pairwise, no relation is between $\TA_\bullet^>(e)$ and $\TA_\bullet^<(e)$, for example. The antipode, if exists, provide
these missing relations,
\begin{equation}
S=\TA^<_\bullet(e)=\RA^>_\bullet(i)=\BA_\bullet^>(e)=\LA_\bullet^<(i)\,.
\end{equation}

The existence of antipode in a Frobenius double algebra depends on if
left transposition (or equivalently, right transposition) w.r.t. $\Fi_B$, i.e., the map $\End\,_BA_B\to\End\,_BA_B$, $X\mapsto X^<$ (or $X\mapsto X^>$), leaves invariant the subalgebra $\TA_V=\End\,_VA\subset\End\,_BA_B$. Similarly, transposition w.r.t. $\Fi_L$ should leave invariant $\RA_H$. This is the content of the next

\begin{lem} \label{lem: S iff}
Let $A$ be a Frobenius double algebra. Then
\begin{itemize}
\item[(a)] we have the following equivalences
\begin{alignat*}{2}
\TA_V^<&\subset \TA_V\quad&\Leftrightarrow\quad \TA_V^>&\subset \TA_V\\
\LA_H^<&\subset \LA_H\quad&\Leftrightarrow\quad \LA_H^>&\subset \LA_H\\
\BA_V^<&\subset \BA_V\quad&\Leftrightarrow\quad \BA_V^>&\subset \BA_V\\
\RA_H^<&\subset \RA_H\quad&\Leftrightarrow\quad \RA_H^>&\subset \RA_H
\end{alignat*}
\item[(b)] $A$ has antipode
\begin{alignat*}{2}
\Leftrightarrow&\quad\TA_V^<\subset \TA_V&\quad\text{and}\quad&\LA_H^<\subset \LA_H\\
\Leftrightarrow&\quad\TA_V^<\subset \TA_V&\quad\text{and}\quad&\RA_H^<\subset \RA_H\\
\Leftrightarrow&\quad\BA_V^<\subset \BA_V&\quad\text{and}\quad&\LA_H^<\subset \LA_H\\
\Leftrightarrow&\quad\BA_V^<\subset \BA_V&\quad\text{and}\quad&\RA_H^<\subset \RA_H
\end{alignat*}
\end{itemize}
\end{lem}
\begin{proof}
\textbf{(a)} $\TA^<_V\subset \TA_V\Rightarrow\exists S\colon A\to A$ s.t. $\TA^<_a=\TA_{S(a)}$, $\forall a\in A$. Then $S=\TA^<_\bullet(e)$ which is invertible by Lemma \ref{lem: corner S}. Thus $\TA^>_a=\TA_{S(S^{-1}(a))}^>=(\TA_{S^{-1}(a)}^<)^>=\TA_{S^{-1}(a)}$ and therefore $\TA_V^>\subset \TA_V$. The backward implication is analogous. The remaining equivalences then follow by symmetry reasons.

\textbf{(b)} If $S$ exists then all transpositions leave invariant their corresponding regular actions. If $\TA_V^<\subset \TA_V$ and $\LA_H^>\subset \LA_H$ then $\TA_a^<=\TA_{S(a)}$ and $\LA_a^<=\LA_{S^{-1}(a)}$ where $S=\TA^<_\bullet(e)$ and $S^{-1}=\LA^>_\bullet(i)$.
Therefore $S$ is an anti algebra automorphism of both $H$ and $V$. Hence $S\Fi_BS^{-1}(a)=S(S^{-1}(a)\ci i)=i\ci a=\Fi_T(a)$,  $S\Fi_LS^{-1}(a)=S(S^{-1}(a)\star i)=i\star a=\Fi_R(a)$ and we have
\begin{align*}
\Fi_T((a\ci a')\star a'')&=S\Fi_B(S^{-1}(a'')\star(S^{-1}(a')\ci S^{-1}(a)))\\
&=S\Fi_B((S^{-1}(a'')\ci a)\star S^{-1}(a'))=\Fi_T(a'\star (S(a)\ci a''))
\end{align*}
One can obtain the analogous relation for $\Fi_R$ similarly.
\end{proof}

After these preparations the following Lemma can be stated without proof.

\begin{lem} \label{lem: S}
Let $S$ be the antipode of a Frobenius double algebra $A$.
\begin{enumerate}
\item $S$ is a double algebra isomorphism $A\iso A^\op_\coop$, i.e.,
\[
S(a\ci a')=S(a')\ci S(a),\ S(e)=e,\ S(a\star a')=S(a')\star S(a),\ S(i)=i.
\]
\item The defining properties in terms of dual bases read as \label{eq: S on qb}
\begin{alignat*}{2}
u_k\oB (v_k\ci S(a))&=(u_k\ci a)\oB v_k&\quad
(x_j\star S(a))\oL y_j&=x_j\oL (y_j\star a)\\
x^j\oR (S(a)\star y^j)&=(a\star x^j)\oR y^j&\quad
(S(a)\ci u^k)\oT v^k&=u^k\oT (a\ci v^k).
\end{alignat*}
\item The restrictions of $S$ to the base ideals are given by
\begin{alignat*}{2}
      S(b)&=\Fi_T\Fi_L(b),\ b\in B,&\quad
      S(l)&=\Fi_R\Fi_T(l),\ l\in L,\\
      S(r)&=\Fi_L\Fi_B(r),\ r\in R,&\quad
      S(t)&=\Fi_B\Fi_R(t),\ t\in T.
\end{alignat*}
\item The following maps are meaningful at least in $\M_k$:
\begin{alignat*}{2}
a\oB a'&\mapsto S(a)\oL a'&\qquad a\oB a'&\mapsto S^{-1}(a')\oR a\\
a\oL a'&\mapsto S^{-1}(a)\oB a'&\qquad a\oL a'&\mapsto S(a')\oT a\\
a\oT a'&\mapsto a\oR S(a')&\qquad a\oT a'&\mapsto a'\oL S^{-1}(a)\\
a\oR a'&\mapsto a\oT S^{-1}(a')&\qquad a\oR a'&\mapsto a'\oB S(a)
\end{alignat*}
\item The antipode relates the dual bases as follows: \label{S rel qb}
\begin{alignat*}{2}
x_j\oL y_j&=S(u_k)\oL v_k&\qquad x^j\oR y^j&=S^{-1}(v_k)\oR u_k\\
u^k\oT v^k&=S(v_k)\oT S(u_k)&\qquad u_k\oB v_k&=S^2(u_k)\oB S^2(v_k).
\end{alignat*}
\end{enumerate}
\end{lem}

For the existence of antipode the following Proposition is useful because it contains
no existential quantifiers, only relations between the structure maps of the Frobenius
double algebra. Unfortunately, it provides only a sufficient condition.

\begin{pro} \label{pro: S iff}
The Frobenius double algebra $A$ has an antipode if
\begin{align}
(u_k\ci a)\star(v_k\ci a')&=\Fi_T(a\star a') \label{eq: 2FT}\\
(x_j\star a)\ci(y_j\star a')&=\Fi_R(a\ci a') \label{eq: 2FR}\\
(a\star x^j)\ci(a'\star y^j)&=\Fi_L(a\ci a') \label{eq: 2FL}\\
(a\ci u^k)\star(a'\ci v^k)&=\Fi_B(a\star a') \label{eq: 2FB}
\end{align}
hold true for all $a,a'\in A$.
If the Galois maps are invertible then
these conditions are also necessary for the existence of $S$.
\end{pro}
\begin{proof}
Using (\ref{eq: 2FB}) one can write
\begin{align*}
\Fi_B(a\star \TA_{a'}^>(a''))&=\Fi_B((a\ci a')\star a'')=
(a\ci a'\ci u^k)\star (a''\ci v^k)\\
&=(a\ci u^k)\star(a''\ci \BA_{a'}^<(v^k))
\end{align*}
Therefore, using also (\ref{eq: 2FT}) in the 3$^{\text{rd}}$ equation
\begin{align*}
\TA_{a'}^>(a'')&=u_i\star\Fi_B(v_i\star \TA_{a'}^>(a''))=u_i\star(v_i\ci u^k)\star
(a''\ci \BA_{a'}^<(v^k))=\\
&=\Fi_T\Fi_R(u^k)\star(a''\ci \BA_{a'}^<(v^k))=a''\ci(\Fi_T\Fi_R(u^k)\star \BA_{a'}^<(v^k))\,,
\end{align*}
hence $\TA_{a'}^>$ is a left $V$-module map, i.e., belongs to $\TA_A$.
Similarly, one can show that (\ref{eq: 2FR}) and (\ref{eq: 2FL}) imply that $\LA^>_A\subset \LA_A$, so by Lemma \ref{lem: S iff}
antipode exists.

Now assume that the conditions of Lemma \ref{lem: G} are satisfied. Since those equations are special cases of equations (\ref{eq: 2FT}), (\ref{eq: 2FL})\dots etc,
sufficiency is obvious. Necessity can be seen also very easily by inserting a product $a\ci a'$ or $a\star a'$ in (\ref{eq: G1}), (\ref{eq: G2}), \dots etc and then using the defining properties of the antipode.
\end{proof}

We remark that if antipode exists then the horizontal multiplication can be expressed as follows
\begin{align}
a\star a'&=a\star u^k\star\Fi_T(v^k\star a')=(a\star u^k)\ci\Fi_R(i\ci(v^k\star a'))
\notag\\
&=(a\star u^k)\ci\Fi_R((S(v^k)\star i)\ci a')\notag\\
&=a\oneT\ci\Fi_R(S(a\twoT)\ci a') \label{eq: clue1}
\end{align}
\begin{align}
a\star a'&=\Fi_T(a\star u^k)\star v^k\star a'=(v^k\star a')\ci\Fi_L(i\ci(a\star u^k))
\notag\\
&=(v^k\star a')\ci\Fi_L((i\star S^{-1}(u^k))\ci a)\notag\\
&={a'}\twoT\ci\Fi_L(S^{-1}({a'}\oneT)\ci a) \label{eq: clue2}
\end{align}
There are two similar expressions that use $\cop_B$ instead. These formulae give a clue to construct the double algebra of a quantum groupoid in Subsection \ref{ss: H}, \ref{ss: WHA} and \ref{ss: HGD}.

\section{Dualities}      \label{s: duals}

\subsection{Duals of almost bialgebroids}

We study dualities between \textit{almost bialgebroids}, the structures found in Proposition \ref{pro: almost QGD}.
If double algebras with antipode are thought of as the analogues of Ocneanu's paragroups then these structures are `paragroups without antipode'.

Let $\bra \A,\B,s,t,\cop,\eps\ket$ be a left almost bialgebroid, that is to say, a left bialgebroid with a possibly non-multiplicative coproduct.  We assume that $\A$ is finitely generated projective both as left and right $\B$-module. Then the two duals of the bimodule $_\B\A_\B$ carry almost right bialgebroid structures in a way the bialgebroids do \cite{KSz}.

The $k$-modules underlying the two duals $\OL{\A}$ and $\OR{\A}$ are
\begin{align*}
\OL{\A}=\Hom(\A_\B,\B_\B)&=
\{\phi\colon\A\to\B\,|\, \phi(t(b)a)=\phi(a)b,\ a\in\A,b\in\B\}\\
\OR{\A}=\Hom(\,_\B\A,\,_\B\B)&=
\{\phi\colon\A\to\B\,|\,\phi(s(b)a)=b\phi(a),\ a\in\A,b\in\B\}
\end{align*}
endowed with algebra structures
\begin{align}
(\phi\psi)(a)&=\psi(s(\phi(a\oneB))\,a\twoB)
\qquad \phi,\psi\in\OL{\A}\\
(\phi\psi)(a)&=\psi(t(\phi(a\twoB))\,a\oneB)
\qquad \phi,\psi\in\OR{\A}
\end{align}
with the unit element being $\eps$ in both cases.
The $\B$-$\B$-bimodule structure of these duals are defined by means of the
source and target homomorphisms
\begin{alignat*}{2}
\OL{s}(b)(a)&=\eps(at(b))  &\qquad \OR{s}(b)(a)&=\eps(a)b \\
\OL{t}(b)(a)&=b\eps(a)  &\qquad \OR{t}(b)(a)&=\eps(as(b))
\end{alignat*}
using the right bialgebroid convention of multiplying with source and target from the right, i.e.,
\begin{equation}
(b\cdot \psi\cdot b')(a)=
\begin{cases}
b\,\psi(a\, t(b'))&\text{ for }\psi\in\OL{\A}\\
\psi(a\, s(b))\,b'&\text{ for }\psi\in\OR{\A}
\end{cases}
\end{equation}
The counits are defined as
\begin{equation*}
\OL{\eps}(\psi)=\psi(1)\qquad\OR{\eps}(\psi)=\psi(1)
\end{equation*}
and are $B$-$B$-bimodule maps in the respective senses.
Moreover, they satisfy
\begin{equation*}
\OR{\eps}(\OR{s}(\OR{\eps}(\phi))\psi)=
\psi(t(\eps(1\twoB)\phi(1))1\oneB)=\psi(t(\phi(1))
1)=\OR{\eps}(\phi\psi)
\end{equation*}
and the analogous one for $\OL{\eps}$, which is one of the right bialgebroid axioms.
Thank to finite projectivity dual comultiplications can be introduced by
\begin{align*}
\psi(aa')&=\psi\oneT(a\,t(\psi\twoT(a')))\qquad \psi\in\OL{\A}\\
\psi(aa')&=\psi\twoT(a\,s(\psi\oneT(a')))\qquad\psi\in\OR{\A}
\end{align*}
which obviously make $\OL{\A}$, respectively $\OR{\A}$, into a comonoid in $_\B\M_\B$. It remains to show the Takeuchi property which for the right dual goes as follows.
\begin{gather*}
\left((\OR{s}(b)\psi\oneT)(a')\cdot\psi\twoT\right)(a)=\left(\psi\oneT(t(b)a')\cdot\psi\twoT\right)(a)=\psi(at(b)a')\\
=\left(\psi\oneT(a')\cdot\psi\twoT\right)(at(b))=\left(\psi\oneT(a')\cdot(\OR{t}(b)\psi\twoT)\right)(a)
\end{gather*}

The properties of these duals are summarized in the next table where we use the
clearer notation $\bra \psi,a\ket$ for $\psi(a)$ if $\psi\in\OL{\A}$ and $\bra a,\psi\ket$ for $\psi(a)$ if $\psi\in\OR{\A}$.
\begin{table}
\begin{tabular}{c|c}
$\bra\under,\under\ket\colon\OL{\A}\x\A\to \B$ &
$\bra\under,\under\ket\colon \A\x\OR{\A}\to \B$\\ \hline
$\bra \phi\psi,a\ket = \bra \psi,\bra\phi,a\oneB\ket\cdot a\twoB\ket$ &
$\bra a,\phi\psi\ket = \bra a\oneB\cdot\bra a\twoB,\phi\ket,\psi\ket$\\
$\bra\psi\OL{s}(b),a\ket = \bra\psi,at(b)\ket$ &
$\bra a,\psi\OR{s}(b)\ket = \bra a,\psi\ket\,b$\\
$\bra\psi\OL{t}(b),a\ket = b\,\bra\psi,a\ket$ &
$\bra a,\psi\OR{t}(b)\ket = \bra as(b),\psi\ket$\\
$\bra\psi,s(b)a\ket = \bra \OL{t}(b)\psi,a\ket$ &
$\bra s(b)a,\psi\ket = b\,\bra a,\psi\ket$\\
$\bra\psi,t(b)a\ket = \bra\psi,a\ket\, b$ &
$\bra t(b)a,\psi\ket = \bra a,\OR{s}(b)\psi\ket$ \\
$\bra\psi, as(b)\ket = \bra\OL{s}(b)\psi,a\ket$ &
$\bra at(b),\psi\ket = \bra a,\OR{t}(b)\psi\ket$\\
$\bra\psi,aa'\ket = \bra\psi\oneT\cdot\bra\psi\twoT,a'\ket,a\ket$ &
$\bra aa',\psi\ket = \bra a,\bra a',\psi\oneT\ket\cdot\psi\twoT\ket$\\
&
\end{tabular}
\caption{The canonical pairings for the two duals of a left `paragroup' $\A$}
\label{tab: pairing}
\end{table}
These pairings can be used also to define duals $\OL{\A'}$ and $\OR{\A'}$ of a right almost bialgebroid $\A'$ so that there will be natural isomorphisms $\OR{\OL{\A}}\cong \A$ and $\OL{\OR{\A}}\cong\A$ for either left or right almost bialgebroids.

The table should make it clear also that $\OR{\A}\cong(\OL{\A_\coop})_\coop$.

Now let $A$ be a Frobenius double algebra. We can apply the above constructions to the left almost bialgebroid
\[
V_B=\bra V,B,\Fi_L,\Fi_R,\cop_B,\Fi_B\ket
\]
found in Proposition \ref{pro: almost QGD}. We will also need the left almost bialgebroid $H_L$. Define
\begin{align}
\OL{\kappa}\colon& H\to \OL{V_B},\qquad h\mapsto\Fi_B(h\star\under)\\
\OR{\kappa}\colon& H\to \OR{V_B},\qquad h\mapsto\Fi_B\Fi_L(\under\ci h)
\end{align}
both of which are $k$-module isomorphisms because $\Fi_B$ is Frobenius. But they preserve more structures,
\begin{align*}
\bra\OL{\kappa}(h)\OL{\kappa}(h'),v\ket&=\Fi_B(h'\star\Fi_B(h\star v\oneB)\star v\twoB)=\Fi_B(h'\star h\star v)\\
&=\bra \OL{\kappa}(h'\star h),v\ket\\
\bra v\ci v',\OR{\kappa}(h)\ket&=\Fi_B\Fi_L(v\ci(\Fi_B\Fi_L(v'\ci h\oneL)\star h\twoL))\\
&=\bra v, \bra v',\OR{\kappa}(h)\oneT\cdot \OR{\kappa}(h)\twoT\ket
\end{align*}
This suggests that the duals of $V_B$ should be closely related to the right almost bialgebroid
\[
H_L^\op=\bra H^\op,L,\Fi_T,\Fi_B,\cop_L,\Fi_L\ket\,.
\]
Preservation of the comultiplication by $\OL{\kappa}$ and the multiplication by $\OR{\kappa}$, however, are not automatic.
\begin{pro}
For any Frobenius double algebra $\bra A,\ci,e,\star,i\ket$ the pairs of maps
$(\OL{\kappa},\Fi_B|_L)$ and $(\OR{\kappa},\Fi_B|_L)$ are isomorphisms of right almost bialgebroids
\[
\begin{CD}
H_L^\op@>\OL{\kappa}>\sim> \OL{V}_B\\
@AAA @AAA\\
L@>\Fi_B|_L>\sim> B
\end{CD}
\qquad
\begin{CD}
H_L^\op@>\OR{\kappa}>\sim> \OR{V}_B\\
@AAA @AAA\\
L@>\Fi_B|_L>\sim> B
\end{CD}
\]
if and only if antipode exists in $A$.
\end{pro}
\begin{proof}
$\OL{\kappa}$ is such an isomorphism iff $\cop_L$ satisfies
\begin{gather*}
\Fi_B(h\star(v\ci v'))=\Fi_B((h\oneL\ci\Fi_L\Fi_B(h\twoL\star v'))\star v)\\
\Leftrightarrow\Fi_L(h\star u_k\ci w))\star v_k=h\oneL\ci\Fi_L\Fi_B(h\twoL\star w)\\
\Leftrightarrow
\Fi_B(h\star(u_k\ci w))\star v_k\oL w'=h\oneL\oL \Fi_B(h\twoL\star w)\star w'\\
\Leftrightarrow
\Fi_B(h\star(u_k\ci u_j))\star v_k\oL v_j=h\oneL\oL h\twoL\\
\Leftrightarrow
\TA^<_{u_j}(h)\oL v_j=h\oneL\oL h\twoL
\end{gather*}
where in each line the quantifier $\forall h,w,w'\in \A$,\dots etc are suppressed. Only the $\Leftarrow$ part of the third $\Leftrightarrow$ needs explanation. Use that both $\under \star w$ and $\Fi_B$ are left $L$ module maps  and then axiom \textbf{A3} to produce $w'$ on the right of $\oL$. Then \textbf{A2} brings $\Fi_T\Fi_L\Fi_B(v_j\star w)$ to the left hand side of $v_k$. Finally use the Nakayama automorphism of $\Fi_B$.
In this way we have proven that the first pair is an isomorphism of right almost bialgebroids iff
\[
\TA^<_A\subset \TA_A\quad\text{and}\quad \TA^<_{u_j}(e)\oL v_j=x_j\oL y_j\,.
\]
For the $\OR{\kappa}$ we obtain that it preserves the almost bialgebroid structures iff it preserves multiplication, i.e., iff
\begin{gather*}
\Fi_B\Fi_L(v\ci(h'\star h))=\Fi_B\Fi_L((v\oneB\star\Fi_B\Fi_L(v\twoB\ci h))\ci h')\\
\Leftrightarrow
\Fi_L(v\ci(x_j\star h))\ci y_j=v\oneB\star\Fi_B\Fi_L(v\twoB\ci h)\\
\Leftrightarrow \dots \dots\\
\Leftrightarrow
\RA^<_A\subset\RA_A\quad\text{and}\quad \RA^<_{x_j}\oB y_j=u_k\oB v_k
\end{gather*}
where we omitted the intermediate steps because the whole argument is the `dual' (vertical $\leftrightarrow$ horizontal) of the previous one.
By Lemma \ref{lem: S iff} (b), second row, we conclude that antipode exists in which case the dual basis relations are automatic (see Lemma \ref{lem: S} (5)).
\end{proof}

If antipode exists in a Frobenius double algebra then there are many dualities between its four almost bialgebroids $V_B$, $H_L$, $V_T$, $H_R$ and their opposites (see
Table \ref{tab: dualities}).
\begin{table}
\begin{tabular}{ccccc}
left `paragroups'& &right `paragroups'&& pairing\\ \hline
$\OR{H^\op_L}\cong V_B$ & and & $H_L^\op\cong\OL{V}_B$ & by &
$\bra h,v\ket_{LB} = \Fi_L\Fi_B(h\star v)$\\
$V_B\cong\OL{H^\op_L}$ & and & $\OR{V}_B\cong H_L^\op$ & by &
$\bra v,h\ket_{BL} = \Fi_B\Fi_L(v\ci h)$\\
$\OR{H}_R\cong V^\op_T$ & and & $H_R\cong \OL{V^\op_T}$ & by &
$\bra h, v\ket_{TR}:=\Fi_T\Fi_R(h\ci v)$\\
$V_T^\op\cong\OL{H}_R$ & and & $\OR{V_T^\op}\cong H_R$ & by &
$\bra v, h\ket_{RT}:=\Fi_R\Fi_T(v\star h)$\\
\hline
&&&&
\end{tabular}
\caption{Dualities inside a double algebra with antipode} \label{tab: dualities}
\end{table}

\subsection{Frobenius integrals}

Left (right) integrals are meaningful in left (right) almost bialgebroids.
Therefore we call an element $i$ of a left almost bialgebroid $\A$ a \textit{left integral} if
$ai=\eps(a)\cdot i\equiv s(\eps(a))i$, $a\in \A$. A \textit{right integral on} $\A$, in turn, is a right integral in $\OR{\A}$, i.e., a $\rho\in\OR{\A}$ satisfying
\begin{equation}
\rho\psi=\rho\OR{s}(\OR{\eps}(\psi))\equiv\rho(\under)\psi(1)\qquad \psi\in\OR{\A}\,.
\end{equation}
Such a $\rho$ is also a right $\B$-module map in the sense of satisfying
\[
\rho(\under s(b))=\rho\OR{t}(b)=\rho\OR{s}(\OR{\eps}(\OR{t}(b)))=\rho(\under)(\OR{t}(b))(1)
=\rho(\under)\,b\,.
\]
One calls $\rho\in\OR{\A}$ a \textit{Frobenius right integral} if it is a right integral on $\A$ and a Frobenius homomorphism for the algebra extension $s\colon\B\to\A$.
Similarly, a left integral $i$ in $\A$, as a functional on $\OR{\A}$, is a $\B$-$\B$ bimodule map because not only the second row in the second column of Table \ref{tab: pairing} holds but
\[
\bra i,\OR{s}(b)\psi\ket=\bra t(b)i,\psi\ket=\bra s(\eps(t(b)))\,i,\psi\ket
=\bra s(b)\,i,\psi\ket=b\,\bra i,\psi\ket
\]
as well. The $i$ is called a \textit{Frobenius left integral} in $\A$ if it is a left integral and a Frobenius homomorphism for the extension $\OR{s}\colon\B\to\OR{\A}$ of algebras. It follows then from standard Frobenius algebra theory that a right integral $\rho$ on $\A$ is a Frobenius homomorphism iff the map
\[
\F\colon \A\to\OR{\A},\qquad a\mapsto \rho(\under a)
\]
is a $k$-module isomorphism. In this case $\F$ is also an isomorphism of $\A$-$\B$-bimodules in the appropriate sense. Note that $\F(a)$ is the analogue of the familiar $a\la\rho$, but "$\rho\ra a$" is not meaningful.

For a Frobenius right integral $\rho$ define $i:=\F^{-1}(\eps)$ which is a left integral because
\begin{align*}
\F(ai)(a')&=\rho(a'ai)=\F(a'a)=\eps(a'a)=\eps(a'\,s(\eps(a)))\\
&=\F(s(\eps(a))\,i)(a')\,,\qquad a,a'\in\A.
\end{align*}
It is called the dual left integral of $\rho$.
If antipode exists, so we are dealing with Hopf algebroids for example, then standard methods show that $i$ is Frobenius. In the general `paragroup' situation this would require more effort which we cannot afford here.

However, all of these concepts become amazingly simple in the double algebraic context.
Left, respectively, right integrals in the almost bialgebroids $V_B$, $V_T$, $H_L$ and $H_R$ have been identified in Lemma \ref{lem: integrals} with the top, bottom, right and left base ideals of $A$. The Frobenius left integrals in $V_B$, for example, are precisely the elements of $T$ that are horizontally invertible, as we show next.
\begin{lem} \label{lem: Frob int}
Let $A$ be a Frobenius double algebra and let $V_B$ and $H_L$ be the underlying left and $V_T$, $H_R$ the underlying right almost bialgebroids. Then
\begin{itemize}
\item the Frobenius left/right integrals in $V_B$/$V_T$ are the (horizontally) invertible elements of $T$/$B$, denoted $T_\star$/$B_\star$,
\item the Frobenius left/right integrals in $H_L$/$H_R$ are the (vertically) invertible elements of $R$/$L$, denoted $R_\ci$/$L_\ci$.
\end{itemize}
If $A$ has an antipode then duality of Frobenius integrals reads as follows.
\begin{itemize}
\item $T_\star$ and $R_\ci$ are in bijection via $r=\Fi_R(t^{-1})$ and $t=\Fi_T(r^{-1})$
\item $B_\star$ and $L_\ci$ are in bijection via $l=\Fi_L(b^{-1})$ and $b=\Fi_B(l^{-1})$
\end{itemize}
\end{lem}
\begin{proof}
Let $t\in V_B$ be a left integral. Then $t\in T$ and $t$ is Frobenius iff the map $\psi\mapsto \bra t,\psi\ket$ is a Frobenius homomorphism $\OR{V}_B\to B$ for the extension $\OR{s}$. Using that the map in Table \ref{tab: dualities} given by $h\mapsto \bra \under,h\ket_{BL}$ is a $k$-module isomorphism $H\to\OR{V}_B$, even though not an isomorphism of almost bialgebroids, we obtain that
$h\mapsto\Fi_T\Fi_L(t\ci h)$ should be a Frobenius homomorphism on the horizontal algebra $H$. Writing $t=i\ci \Fi_L(t)$ this is equivalent to that
\[
h\mapsto \Fi_T(\Fi_L(t)\ci h)=\Fi_T(\Fi_B\Fi_L(t)\star h)
\]
be a Frobenius homomorphism for the subalgebra $T\subset H$. But $\Fi_T$ is Frobenius by assumption so $t$ is Frobenius iff $\Fi_B\Fi_L(t)$ is invertible in $\Cnt_H(T)$. Due to Lemma \ref{lem: inv} this is equivalent to that $\Fi_B\Fi_L(t)$ is invertible in $B$, i.e., $t$ is invertible in $T$, i.e., $t\in T_\star$.
Similar arguments work for $B$, $R$ and $L$.

For duality of Frobenius integrals consider again the duality $V_B\leftrightarrow H_L^\op$ provided by $\bra \under,\under\ket_{BL}$ and by the isomorphism $\Fi_B|_L$.
The Frobenius left integrals in $V_B$ are the elements of $T_\star$ and the Frobenius left integrals in $H_L$ are the elements of $R_\ci$. For $t\in T_\star$ and $r\in R$
\[
\bra t,r\star h\ket_{BL}=\Fi_B\Fi_L(h)\,\ \forall h\in H\quad\Leftrightarrow\quad
t\star S^{-1}(r)=e\,.
\]
Using that $S^{-1}|_R=\Fi_L\Fi_T|_R$ this is equivalent to $t\star \Fi_T(r)=\Fi_T(e)=i$, i.e. that $t^{-1}=\Fi_T(r)$. But $\Fi_T|_R$ is an algebra isomorphism, so $r\in R_\ci$. Now let $r\in R_\ci$ and $t\in T$ then
\[
\bra v\ci t,r\ket_{BL}=\Fi_B(v),\ \forall v\in V\quad\Leftrightarrow\quad
t\ci r=i\,.
\]
Therefore the dual integral of $r$ is the $t\in T$ satisfying $\Fi_R(t)\ci r=\Fi_R(i)=e$, i.e., $t=\Fi_T(r^{-1})$. This is just the inverse of the previous
construction of $r$ from $t$. The remaining dualities are left to the reader.
\end{proof}

\section{Distributive Frobenius double algebras} \label{s: D}

The compatibility conditions between the vertical and horizontal multiplications in  Frobenius double algebras are too weak to ensure multiplicativity of comultiplication or the existence of antipode. So we need some further assumption in order to obtain
Hopf algebroids. No doubt, the most natural compatibility between two monoid structures is distributivity. What we apply here, however, involves
also the comultiplications. It should be understood, therefore, as distributivity between two Frobenius structures.

\begin{defi} \label{def: D}
A Frobenius double algebra $A$ is called distributive if for all $a,a',a''\in A$
\begin{align}
a\ci(a'\star a'')&=(a\oneB\ci a')\star(a\twoB\ci a'') \label{eq: DB}\\
a\star(a'\ci a'')&=(a\oneL\star a')\ci(a\twoL\star a'') \label{eq: DL}\\
(a'\star a'')\ci a&=(a'\ci a\oneT)\star(a''\ci a\twoT) \label{eq: DT}\\
(a'\ci a'')\star a&=(a'\star a\oneR)\ci(a''\star a\twoR) \label{eq: DR}
\end{align}
\end{defi}
Inserting $a=i$ in the first and third and $a=e$ in the second and fourth distributivity law we obtain equations (\ref{eq: 2FT}), (\ref{eq: 2FB}), (\ref{eq: 2FR}) and (\ref{eq: 2FL}), respectively.
\begin{cor}
In a distributive double algebra antipode exists and the Galois maps $\Gamma_{XY}$ are invertble with inverse $\Gamma_{YX}$.
\end{cor}

The distributive laws can also be interpreted as module algebra properties. For example, (\ref{eq: DB}) means that $\LA$ is a left module algebra action of the algebra $V$, with comultiplication $\Delta_B$, on the algebra $H$.
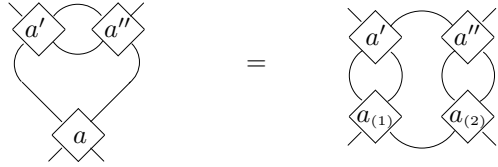
\begin{equation*}
\parbox[c]{1.3in}{
\begin{picture}(70,80)
\put(15,43){\rhomb{\!\!a'}}
\put(45,43){\rhomb{\!\!a''}}
\put(43,55){\arc{23}{0.6}{2.5}}
\put(43,51){\arc{23}{3.8}{5.6}}
\put(30,3){\rhomb{a}}
\put(56.5,40.5){\arc{20}{5.5}{7.2}}
\put(37,18){\line(-1,1){15}}
\put(29,40){\arc{20}{2.4}{4}}
\put(48,18){\line(1,1){15}}
\put(37,8){\line(-1,-1){5}}
\put(48,8){\line(1,-1){5}}
\put(22,58){\line(-1,1){5}}
\put(63,58){\line(1,1){5}}
\end{picture}
}
\quad =\quad
\parbox[c]{1.6in}{
\begin{picture}(100,80)
\put(46,45){\arc{30}{3.8}{5.6}}
\put(15,40){\rhomb{\!a'}}
\put(30,35.5){\arc{23}{2.2}{4.1}}
\put(27,35.5){\arc{23}{5.4}{7.3}}
\put(15,10){\rhomb{\!\!\!\raise1pt\hbox{$a\scriptstyle\oneB$}}}
\put(50,40){\rhomb{\!\!a''}}
\put(65,35.5){\arc{23}{2.2}{4.1}}
\put(62,35.5){\arc{23}{5.4}{7.3}}
\put(50,10){\rhomb{\!\!\!\raise1pt\hbox{$a\scriptstyle\twoB$}}}
\put(46,23){\arc{30}{0.6}{2.5}}
\put(23,56){\line(-1,1){5}}
\put(23,14){\line(-1,-1){5}}
\put(69,56){\line(1,1){5}}
\put(69,14){\line(1,-1){5}}
\end{picture}
}
\end{equation*}

\begin{pro} \label{pro: D}
A Frobenius double algebra $\bra A,\ci,e,\star,i\ket$ is distributive if and only if
antipode exists in $A$ in the sense of Definition \ref{def: S} and the comultiplications are multiplicative, i.e.,
\begin{alignat*}{2}
\Delta_B\colon&V\to V\ex{B}V&\qquad\Delta_R\colon&H\to H\ex{R} H\\
\Delta_L\colon&H\to H\ex{L}H&\qquad\Delta_T\colon&V\to V\ex{T}V
\end{alignat*}
are algebra homomorphisms.
\end{pro}
\begin{proof}
Recall Proposition \ref{pro: almost QGD} that the $\Delta$'s obey the Takeuchi property, hence they are $k$-module maps of the indicated type. Due to the previous
Corollary antipode exists in distributive Frobenius double algebras. Therefore we only have to show that in the presence of antipode distributivity is equivalent to multiplicativity of the $\Delta$'s. Consider $\Delta_B$. Using that $\Fi_B$ is Frobenius the $\Delta_B$ is multiplicative iff
\[
(a'\ci a'')\oneB\star\Fi_B((a'\ci a'')\twoB\star a)=
(a'\oneB\ci a''\oneB)\star\Fi_B((a'\twoB\ci a''\twoB)\star a)
\]
holds for all $a,a',a''\in A$.
Inserting the definition of $\Delta_B$ and using the dual basis property this is equivalent to the equation
\[
(a'\ci a'')\star a=
((a'\star u_j)\ci(a''\star u_k))\star\Fi_B((v_j\ci v_k)\star a)
\]
The LHS is the same as the LHS of (\ref{eq: DR}), so we can concentrate on the RHS.
By means of the antipode we can transpose $v_k$ to the right therefore the RHS can be written as
\begin{gather*}
(a'\star u_j\star\Fi_B((v_j\ci v_k)\star a))\ci(a''\star u_k)=\\
(a'\star u_j\star\Fi_B(v_j\star(a\ci S^{-1}(v_k))))\ci(a''\star u_k)=\\
(a'\star (a\ci S^{-1}(v_k)))\ci(a''\star u_k)=\\
(a'\star (a\ci x^k))\ci(a''\star y^k)
\end{gather*}
where in the last step Lemma \ref{lem: S} (\ref{S rel qb}) has been used. Taking into account the definition of $\Delta_R$, given in Section \ref{s: F2}, this is precisely the RHS of (\ref{eq: DR}). Arguing with the opposite horizontal structure we could show that multiplicativity of $\Delta_B$ is equivalent also to (\ref{eq: DL}), provided the antipode exists. Passing to the opposite vertical structure all of these are equivalent to that $\Delta_T$ is multiplicative. Similarly, under the existence of antipode (\ref{eq: DB}) $\Leftrightarrow$ (\ref{eq: DT}) $\Leftrightarrow$ $\Delta_L$ is multiplicative $\Leftrightarrow$ $\Delta_R$ is multiplicative.
\end{proof}

\begin{thm} \label{thm: main}
Let $\bra A,\ci,e,\star,i\ket$ be a distributive Frobenius double algebra. Then
\begin{itemize}
\item $V$ is a Hopf algebroid \cite{HGD} with underlying left and right bialgebroids
\[
\bra V,B,\Fi_L|_B,\Fi_R|_B,\Delta_B,\Fi_B\ket \quad\text{and}\quad
\bra V,T,\Fi_R|_T,\Fi_L|_T,\Delta_T,\Fi_T\ket
\]
respectively, such that $i$ is a two-sided Frobenius integral in $V$ and the antipode is the double algebraic antipode of $A$.
\item $H$ is a Hopf algebroid with underlying left and right bialgebroids
\[
\bra H,L,\Fi_B|_L,\Fi_T|_L,\Delta_L,\Fi_L\ket \quad\text{and}\quad
\bra H,R,\Fi_T|_R,\Fi_B|_R,\Delta_R,\Fi_R\ket
\]
respectively, such that $e$ is a two-sided Frobenius integral in $H$ and the antipode is the inverse of the double algebraic antipode of $A$.
\item The vertical Hopf algebroid $V$ and the horizontal Hopf algebroid $H$ given above are in duality w.r.t any one of the following pairings:
\begin{alignat*}{2}
\bra h,v \ket_{LB}&:=\Fi_L\Fi_B(h\star v)&\qquad
\bra v, h \ket_{BL}&:=\Fi_B\Fi_L(v\ci h )\\
\bra v,h\ket_{RT}&:=\Fi_R\Fi_T(v\star h)&\qquad
\bra h,v \ket_{TR}&=\Fi_T\Fi_R(h\ci v)
\end{alignat*}
\end{itemize}
Conversely, every Hopf algebroid possessing a two-sided Frobenius integral is the vertical (or horizontal) Hopf algebroid of a distributive Frobenius double algebra $A$.
\end{thm}

\begin{proof}
By Proposition \ref{pro: D} and Proposition \ref{pro: almost QGD} the four 6-tuples of the Theorem are all bialgebroids and there exists a (double algebraic) antipode $S$ on $A$. At first we will show that the axioms of a Hopf algebroid antipode given in \cite[Definition 4.1]{HGD} are satisfied by $S$, for both $V$ and $H$. This definition involves only the left bialgebroid and $S$.
Consider the left bialgebroid $V_B$. The first axiom claims that the source
map $s_L=\Fi_L|_B$  and the target map $t_L=\Fi_R|_B$ are related by $St_l=s_L$. But
this is obvious in the double algebra because $S|_R=\Fi_L\Fi_B$.
The remaining two axioms are less trivial but also not difficult calculations within the double algebra:
\begin{gather*}
S(a\oneB)\oneB\ci a\twoB\oB S(a\oneB)\twoB=(S(a\star u_j)\star u_k)\ci v_j\oB v_k\\
=(x_j\star S(a)\star u_k)\ci y_j\oB v_k=(x_j\star u_k)\ci(y_j\star a)\oB v_k\\
=u_k\ci(y_j\star a)\oB v_k\star x_j=u_k\oB(v_k\ci S(y_j\star a))\star x_j\\
=u_k\oB(v_k\ci(S(a)\star u^j))\star v^j=u_k\oB(v_k\ci u^j)\star v^j\star S(a)\\
=u_k\oB\Fi_B\Fi_L(v_k)\star S(a)=e\oB S(a)
\end{gather*}
\begin{gather*}
S^{-1}(a\twoB)\oneB\oB S^{-1}(a\twoB)\twoB\ci a\oneB=S^{-1}(v_j\star a)\star u_k\oB v_k\ci u_j\\
=S^{-1}(a)\star x^j\star u_k\oB v_k\ci y^j=x^j\star u_k\oB v_k\ci(a\star y^j)\\
=x^j\star(u_k\ci S^{-1}(a\star y^j))\oB v_k=u^j\star(u_k\ci(v^j\star S^{-1}(a)))\oB v_k\\
=S^{-1}(a)\star u^j\star(u_k\ci v^j)\oB v_k=S^{-1}(a)\star\Fi_B\Fi_R(u_k)\oB v_k\\
=S^{-1}(a)\oB e
\end{gather*}
This proves that $V_B$ is a Hopf algebroid with antipode $S$. The dual calculation, in which vertical and horizontal are interchanged, proves that $H_L$ is a Hopf algebroid with antipode $S^{-1}$.
But then it is an easy exercise for the reader to check that the right bialgebroids $V_T$ and $H_R$ are images under $S$ of the left bialgebroid structures, so that $(V_B,S,V_T)$ and $(H_L,S^{-1},H_R)$ are Hopf algebroids \cite{HGD} in `symmetrized form'.

$i$ being an element of both $B_\star$ and $T_\star$ it is a 2-sided Frobenius integral of the Hopf algebroid $V$ by Lemma \ref{lem: Frob int}.

As for the dualities we refer to Section \ref{s: duals} where the four pairings were shown to provide dualities between the underlying almost bialgebroids in Table \ref{tab: dualities}. Therefore the 4 underlying bialgebroids of $A$ are in duality in the same sense: The pairings satisfy the relations listed in Table \ref{tab: pairing}.
We remark also that the antipode relates these pairings as
\begin{align*}
S(\bra h,v\ket_{LB})&=\bra S(v),S(h)\ket_{RT}\\
S(\bra v,h\ket_{BL})&=\bra S(h),S(v)\ket_{TR}
\end{align*}
which are simple consequences of that the antipode is a double algebra antiautomorphism.

The proof of the assertion that every Hopf algebroid with a two-sided Frobenius integral is the vertical Hopf algebroid of a double algebra, is altogether shifted to
Subsection \ref{ss: HGD}.
\end{proof}

\section{Examples}

\subsection{Commutative algebras}

For a $k$-algebra $A$ with multiplication $\bra x,y\ket\mapsto xy$ and unit element $1$ define $x\ci y:=xy$, $x\star y:=xy$ and $e:=1=:i$. Then $\bra A,\ci,e,\star,i\ket$ is a double algebra precisely if the original algebra is commutative.

\subsection{Full matrix algebras}

Let $A=\mathtt{M}_n(k)$ with matrix units $\{e_{jk}\}$ and define the
associative operations
\begin{align*}
e_{jk}\ci e_{lm}&=\delta_{k,l}e_{jm}\\
e_{jk}\star e_{lm}&=\delta_{j,l}\delta_{k,m} e_{jk}
\end{align*}
where $\delta$ denotes the Kronecker symbol. They have units $e=\sum_j e_{jj}$ and $i=\sum_{jk} e_{jk}$, respectively. Then $A$ becomes a double algebra with base homomorphisms
\begin{align*}
\Fi_L(e_{jk})&=\Fi_R(e_{jk})=\delta_{j,k}e_{jk}\\
\Fi_B(e_{jk})&=\sum_l e_{jl}\\
\Fi_T(e_{jk})&=\sum_l e_{lk}
\end{align*}
and with antipode $S(e_{jk})=e_{kj}$. This is a special case of the next groupoid example which, in turn, is a special case of weak Hopf algebras.

\subsection{Finite groupoids}

Let $k$ be a field and $s,t\colon G\rightrightarrows O$ be a finite groupoid. Defining
$A=kG$, the $k$-vector space with basis $G$, and
\begin{equation*}
\begin{aligned}
g\ci g'&=gg'\\
g\star g'&=\delta_{g,g'}g
\end{aligned}
\quad\qquad
\begin{aligned}
e&=\sum_{x\in O}x\\
i&=\sum_{g\in G}g
\end{aligned}
\quad\qquad S(g)=g^{-1}
\end{equation*}
the $\bra A,\ci,e,\star,i\ket$ is a double algebra with antipode. The base homomorphisms are
\begin{align*}
\Fi_L(g)&=\Fi_R(g)=\delta_{g\in O} g\\
\Fi_B(g)&=\sum_{g'\in G,\ t(g')=t(g)}\ g'\\
\Fi_T(g)&=\sum_{g'\in G,\ s(g')=s(g)}\ g'
\end{align*}
Since finite groupoid algebras are weak Hopf algebras, this is a special case of Subsection \ref{ss: WHA}, hence a distributive Frobenius double algebra.

\subsection{Finite double categories}

In a recent paper N. Andruskiewitsch and S. Natale \cite{AN} have shown that finite vacant double groupoids have a natural weak Hopf algebra structure. It is natural to expect that the notion of double algebra allows even more double categories.

Let $\D$ be a double category with horizontal 1-cells $\Ha$, vertical 1-cells $\Ve$, 0-cells $\Oo$, horizontal composition $\star$ and vertical composition $\ci$ such that the set $\D$ of 2-cells is finite. For a commutative ring $k$ let $A:=k\D$ be the free $k$-module generated by the set of 2-cells. Then one can extend the compositions $\star$ and $\ci$ to be $k$-linear associative multiplications on $A$ by postulating the multiplication of uncomposable cells to be zero. The horizontal multiplication $\star$ has unit $i=\sum_{v\in\Ve}v$, the sum of vertical 1-cells, and the vertical multiplication $\ci$ has unit $e=\sum_{h\in\Ha}h$ the sum of horizontal 1-cells.
(We consider the 0-cells and 1-cells as subsets of $\D$.)
Denoting by $\beta,\tau\colon\D\to\Ha$ and $\lambda,\rho\colon\D\to\Ve$ the source and target maps of the vertical, respectively, horizontal categories we obtain the following boundary homomorphisms, evaluated on $c\in\D\subset A$:
\begin{align*}
\Fi_B(c)&=c\ci i=\begin{cases}\sum_{v\in\Ve,\ \beta(v)=\tau(c)}c\ci v&
                               \text{if $\tau(c)\in\Oo$}\\
                              \ \ \ 0&\text{otherwise}\end{cases}\\
\Fi_T(c)&=i\ci c=\begin{cases}\sum_{v\in\Ve,\ \tau(v)=\beta(c)}v\ci c&
                               \text{if $\beta(c)\in\Oo$}\\
                              \ \ \ 0&\text{otherwise}\end{cases}\\
\Fi_L(c)&=c\star e=\begin{cases}\sum_{h\in\Ha,\ \lambda(h)=\rho(c)}c\star a&
                                \text{if $\rho(c)\in\Oo$}\\
                                \ \ \ 0&\text{otherwise}\end{cases}\\
\Fi_R(c)&=e\star c=\begin{cases} \sum_{h\in\Ha,\ \rho(v)=\lambda(c)}a\star c&
                              \text{if $\lambda(c)\in \Oo$}\\
                              \ \ \ 0& \text{otherwise}\end{cases}
\end{align*}
These $\Fi(c)$'s are like kites with fixed "triangular" head $c$ and all possible 1-cells as tails and moving downward, upward, left and right, respectively.
Taking into account the $D_4$ symmetry of the axioms of both double categories and
double algebras, in order to see that $A$ is a double algebra it suffices to check axiom \textbf{A1}.
\begin{align*}
\Fi_L(c)\ci d&=
\begin{cases}
\sum_h(c\star h)\ci d& \text{where the sum is over $h\in\Ha$ s.t. $\beta(d)=\tau(c)\star h$}\\
\ \ \ 0 &\text{if such $h$ does not exist}
\end{cases}\\
\Fi_B\Fi_L(c)\star d&=
\begin{cases}
\sum_h((c\star h)\ci v)\star d& \text{where $v=\lambda(d)$ and the sum runs over the}\\
&\text{$h\in\Ha$ s.t. $\tau(c)\star h=\lambda\tau(c)=\beta(v)$}\\
\ \ \ 0& \text{if such $h$ does not exist}
\end{cases}
\end{align*}
We see that these expressions can be equal only if all horizontal 1-cells have right inverses. If this is satisfied then both expessions are equal to
\[
\begin{cases}
(c\star\tau(c)^{-1})\ci d&\text{if $\rho(c)\in\Oo$ and $\beta\lambda(d)=\lambda\tau(c)$}\\
\ \ \ 0&\text{otherwise}
\end{cases}
\]
Applying the symmetry operations we obtain that $A$ is a double algebra precisely if the categories $\Ha$ and $\Ve$ are groupoids. This case, of course, contains the vacant double groupoids but also every double groupoids. The question, however, when $A$ is Frobenius and distributive requires further investigations.

\subsection{Frobenius extensions} \label{ss: F ext}
This is an example of a double algebra which can be neither distributive nor Frobenius but has antipode.

Let $N\subset M$ be a Frobenius extension of $k$-algebras with Frobenius homomorphism $\psi\colon\,_NM_N\to\,_NN_N$ and dual basis $\sum_i e_i\o f_i$. Define $A$ to be the center of $M\o_N M$ considered as an $N$-$N$ bimodule, i.e., $A:=(M\o_N M)^N$. General elements of $A$ are denoted by $a=a_1\o a_2$, with the summation over a finite index set suppressed. One can introduce two algebra structures on $A$ as follows.
\begin{alignat*}{2}
a\ci a'&:=a_1a'_1\o_N a'_2a_2&\qquad e&:=1\o_N 1\\
a\star a'&:=a_1\psi(a_2a'_1)\o_N a'_2&\qquad i&:=e_i\o_N f_i
\end{alignat*}
We claim that this structure on $A$ is a double algebra. At first we compute the base
homomorphisms:
\begin{align*}
\Fi_L(a)&=a\star e=a_1\psi(a_2)\o_N 1\\
\Fi_R(a)&=e\star a=1\o_N \psi(a_1)a_2\\
\Fi_B(a)&=a\ci i=a_1e_i\o_N f_ia_2\\
\Fi_T(a)&=i\ci a=e_ia_1\o_N a_2f_i
\end{align*}
It is now easy to verify that \textbf{A1}-\textbf{A8} are satisfied. We write down
some of them:
\begin{align*}
\Fi_T\Fi_L(a)\star a'&=e_ia_1\psi(a_2)\psi(f_ia'_1)\o a'_2=a'_1a_1\psi(a_2)\o a'_2\\
&=a'\ci\Fi_L(a)\\
\Fi_B\Fi_L(a)\star a'&=a_1\psi(a_2)a'_1\o a'_2=\Fi_L(a)\ci a'\\
\Fi_R\Fi_B(a)\ci a'&=a'_1\o a'_2\psi(a_1e_i)f_ia_2=a'_1\o a'_2\,a_1a_2\\
&=a'_1\psi(a'_2a_1e_i)\o f_ia_2=a'\star\Fi_B(a)\,.
\end{align*}

The antipode of this double algebra and its inverse are
\begin{align*}
S(a)&=\psi(e_k a_1)a_2\o_N f_k\\
S^{-1}(a)&=e_k\o_N a_1\psi(a_2 f_k)\,.
\end{align*}
For example, (\ref{eq: SB}) can be proven by
\begin{align*}
\Fi_B(a'\star(a''\ci a))&=a'_1\psi(a'_2a''_1a_1)e_k\o_N f_ka_2a''_2\\
&=a'_1\tilde a_1\psi(\tilde a_2a'_2a''_1)e_k\o_N f_ka'_2\\
&=\Fi_B((a'\ci S(a))\star a'')
\end{align*}
where $\tilde a_1\o \tilde a_2:=S(a)$ and we used the fact that $\tilde a_1\psi(\tilde a_2 x)=\psi(xa_1)a_2$ holds for all $x\in M$.

The double algebra with antipode (or paragroup?) we have found can be related to the Jones tower $N\subset M\subset M_2\subset M_3$ as follows.
The 2-step relative commutants $\Cnt_{M_2}(N)=\End\,_NM_N$ and $\Cnt_{M_3}(M)=\End\,_M(M\o_N M)_M$ are algebras for composition of endomorphisms.
There are two $k$-module isomorphisms from $A=(M\o_N M)^N$,
\begin{align*}
\pi\colon& A\to\End\,_NM_N\,,\qquad \pi(a)(m)=a_1\psi(a_2m)\\
\vartheta\colon&A\to\End\,_M(M\o_N M)_M\,,\qquad\vartheta(a)(m\o_N m')=ma_1\o_N a_2m'
\end{align*}
the first of which is requiring the Frobenius structure. Pulling back the algebra structures via $\pi$ and $\vartheta$ to $A$ one obtains the horizontal multiplication $\star $ and the opposite of the vertical multiplication $\ci$. In this way both of the 2-step relative commutants are double algebras and are in duality position.

\subsection{Depth 2 Frobenius extensions}

If the Frobenius extension $N\subset M$ discussed in Subsection \ref{ss: F ext} is of depth 2 \cite{KSz} then the double algebra $A$ constructed above is a distributive Frobenius double algebra. Although this is a consequence of (the converse part of) Theorem \ref{thm: main} and of the results of \cite{FrobD2} a direct proof is desirable and follows below.

Recall \cite{KSz} that $N\subset M$ being of depth 2 is equivalent to the existence
of $b_j,c_j\in(M\o_N M)^N$ and $\beta_j,\gamma_j\in\End\,_NM_N$ such that
\begin{align}
b_j^1\o_N b_j^2\beta_j(m)m'&=m\o_N m'\\
m\gamma_j(m')c_j^1\o_N c_j^2&=m\o_N m' \label{eq: right D2}
\end{align}
for all $m\o_N m'\in M\o_N M$.
In case of $N\subset M$ is also Frobenius the D2 bases are related by a Frobenius system. Let $\psi\colon M\to N$ be a Frobenius homomorphism with dual basis $e_k\o_N f_k$ then a right D2 basis can be obtained from a left D2 basis by
\begin{equation*}
\gamma_j(\under)=\psi(\under b_j^1)b_j^2\,,\qquad
c_j^1\o_N c_j^2=\beta_j(e_k)\o_N f_k\,.
\end{equation*}
The presence of such a D2 basis causes the double algebra structure introduced in Subsection \ref{ss: F ext} on $A=(M\o_N M)^N$ to be Frobenius and distibutive as we will show now.

$\Fi_B$ is Frobenius. As a matter of fact, let $u_j:=b_j$ and $v_j:=c_j$. Then
\begin{align*}
\Fi_B(a\star u_j)\star v_j&=a_1\psi(a_2b_j^1)b_j^2c_j^1\o_N c_j^2
=a_1\gamma_j(a_2)c_j^1\o_N c_j^2=a\\
u_j\star\Fi_B(v_j\star a)&=b_j^1\o_N\psi(b_j^2c_j^1\psi(c_j^2a_1)e_k)f_ka_2\\
&=b_j^1\o_N b_j^2\beta_j(e_k)\psi(f_k a_1)a_2=b_j^1\o_N b_j^2\beta_j(a_1)a_2=a\,.
\end{align*}

$\Fi_L$ is Frobenius. As a matter of fact, let
\begin{align*}
x_j&:=S(u_j)=\gamma_j(e_k)\o_N f_k\\
y_j&:=v_j=\beta_j(e_k)\o_N f_k\,.
\end{align*}
Then
\begin{align*}
\Fi_L(a\ci x_j)\ci y_j&=a_1\gamma_j(e_k)\psi(f_ka_2)\beta_j(e_l)\o_N f_l
=a_1\gamma_j(a_2)c_j^1\o_N c_j^2=a\\
x_j\ci\Fi_L(y_j\ci a)&=\gamma_j(e_k)\beta_j(e_l)a_1\psi(a_2f_l)\o_N f_k\\
&=\gamma_j(e_k)c_j^1a_1\psi(a_2c_j^2)\o_N f_k=a_1\psi(a_2e_k)\o_N f_k=a\,.
\end{align*}

Since antipode exists, the $\Fi_T=S\Fi_BS^{-1}$ and $\Fi_R=S\Fi_LS^{-1}$ are also Frobenius. The corresponding comultiplications are
\begin{align*}
\cop_B(a)&=a\star u_j\oB v_j=(a_1\o_N\gamma_j(a_2))\oB(c_j^1\o_N c_j^2)\\
\cop_L(a)&=a\ci x_j\oL y_j=(a_1\gamma_j(a_2e_k)\o_N f_k)\oL(c_j^1\o_N c_j^2)\\
\cop_T(a)&=a\star u^j\oT v^j=(a_1\o_N \beta_j(a_2))\oT(e_k\o_N b_j^1\psi(b_j^2f_k))\\
\cop_R(a)&=a\ci x^j\oR y^j=(e_k\o_N\beta_j(f_ka_1)a_2)\oR(b_j^1\o_N b_j^2)\,.
\end{align*}

Distributivity. Again by the existence of $S$ it suffices to prove (\ref{eq: DB}) and (\ref{eq: DL}).
\begin{align*}
a\ci(a'\star a'')&=a_1a'_1\psi(a'_2a''_1)\o_N a''_2a_2\\
(a\oneB\ci a')\star(a\twoB\ci a'')&=a_1a'_1\psi(a'_2\gamma_j(a_2)c_j^1a''_1)\o_N a''_2c_j^2
\end{align*}
which are indeed equal due to (\ref{eq: right D2}). Left distributivity demands
\begin{align*}
(a\oneL\star a')\ci(a\twoL\star a'')&=(a_1\gamma_j(a_2e_k)\psi(f_ka'_1)\o_N a'_2)
\ci(c_j^1\psi(c_j^2a''_1)\o_N a''_2)\\
&=a_1\gamma_j(a_2a'_1)c_j^1\psi(c_j^2a''_1)\o_N a''_2a'_2\\
&=a_1\psi(a_2a'_1a''_1)\o_N a''_2a'_2
\end{align*}
to be equal to $a\star(a'\ci a'')$ which is clear.

\subsection{Endomorphism monoids of Frobenius objects}

Underlying of the double algebra is the category $\M_k$ of $k$-modules. Replacing it with any symmetric monoidal closed category $\bra\V,\o,I\ket$ we obtain the notion \textbf{double monoids in $\V$}.

Let $\bra f,\mu,\eta,\gamma,\pi\ket$ be a Frobenius algebra in the monoidal category
$\bra\C,\bo,U\ket$. That is to say,
\begin{enumerate}
\item $\bra f,\mu,\eta\ket$ is a monoid,
\item $\bra f,\gamma,\pi\ket$ is a comonoid and
\item the Frobenius properties hold:
\begin{align*}
(\mu\bo f)\ci\asso_{f,f,f}\ci(f\bo\gamma)&=\gamma\ci\mu\\
(f\bo\mu)\ci\asso^{-1}_{f,f,f}\ci(\gamma\bo f)&=\gamma\ci\mu
\end{align*}
\end{enumerate}
The endomorphism monoid $A=\End f$ is a monoid (in $\Set$) with multiplication given by composition $\ci$ and with unit given by the identity arrow $f:f\to f$.
But there is another monoid structure on $A$ given by the convolution product
$a\star b:=\mu\ci(a\bo b)\ci\gamma$, $a,b\in A$, which has unit $\iota:=\eta\ci\pi$.
It is easy to see that the two monoid structures obey the axioms of a double monoid in $\Set$. If $\C$ is a $k$-linear monoidal category then $A$ is a double algebra over $k$. The more general situation of a depth 2 Frobenius arrow in a bicategory has been studied in \cite{FrobD2}.

\subsection{Takeuchi products}

Here we prove that for $A$ any Frobenius double algebra with antipode the Takeuchi
product $A\x_B A$ is a double algebra.

Let $a_1\oB a_2$ denote a general element on $A\oB A$, not just rank 1 tensors. Then the $\x_B$-product
\[
A\x_B A:=\{a_1\oB a_2\in A\oB A\,|\, a_1\oB a_2\ci\Fi_L(b)=
a_1\ci\Fi_R(b)\oB a_2,\ b\in B\,\}
\]
is a ring with respect to termwise vertical multiplication. If $a_1\oB a_2$ belongs to $A\x_B A$ we use the notation $a_1\x_B a_2$. So we define
\[
(a_1\x_B a'_1)\ci(a_2\x_Ba'_2):=(a_1\ci a_2)\x_B(a'_1\ci a'_2)
\]
which clearly makes $A\x_B A$ into a $k$-algebra with unit $e\oB e$. Indeed,
\[
e\oB e\ci\Fi_L(b)=e\oB b\star e=e\star b\oB e=e\ci\Fi_R(b)\oB e
\]
so $e\oB e\in A\x_B A$. But there exists another multiplication on the whole of $A\oB A$,
\[
(a_1\oB a_2)\star(a'_1\oB a'_2):=a_1\star\Fi_B(a_2\star a'_1)\oB a'_2
\]
with unit $u_k\oB v_k$. In this horizontal algebra $A\x_B A$ is a subalgebra. As a matter of fact, using the Nakayama automorphism of $\Fi_B$
\begin{gather*}
a_1\star\Fi_B(a_2\star a'_1)\oB a'_2\ci\Fi_L(b)=
a_1\star\Fi_B(a_2\star a'_1\star\Fi_T\Fi_R(b))\oB a'_2=\\
a_1\star\Fi_B(\Fi_T\Fi_L(b)\star a_2\star a'_1)\oB a'_2=
a_1\ci \Fi_R(b)\oB\Fi_B(a_2\star a'_1)\star a'_2
\end{gather*}
so $A\x_B A$ is closed under $\star$ and it contains the unit because by (\ref{eq: Takeuchi prop}) the $u_k\oB v_k=\Delta(i)\in A\x_B A$. Now we are going to show that these two algebra structures on $A\x_B A$ obey the axioms of double algebras. At first we compute the base homomorphisms which we denote by $\beta$.
\begin{align*}
\beta_L(a_1\x a_2)&:=(a_1\x a_2)\star (e\x e)=a_1\star\Fi_B\Fi_L(a_2)\x e\\
\beta_R(a_1\x a_2)&:=(e\x e)\star(a_1\x a_2)=e\x\Fi_B\Fi_R(a_1)\star a_2\\
\beta_B(a_1\x a_2)&:=(a_1\x a_2)\ci(u_k\x v_k)=(a_1\ci u_k)\x(a_2\ci v_k)\\
\beta_T(a_1\x a_2)&:=(u_k\x v_k)\ci(a_1\x a_2)=(u_k\ci a_1)\x(v_k\ci a_2)
\end{align*}

Replacing $A$ with $A_\coop$ the structure of $A\x_B A$ changes as follows.
There is an isomorphism $A\oB A\iso A_\coop\oB A_\coop$, $a\o a'\mapsto a'\o a$ under which $A\x_B A$ is mapped to $A_\coop\x_B A_\coop$ (because $\Fi_L$ and $\Fi_R$ are interchanged), the vertical multiplication is invariant and the horizontal multiplication changes to its opposite. Therefore axioms A4, A5, A6 and A7 for $A\x_B A$ become the axioms A3, A2, A1 and A8, respectively, for $A_\coop\x_B A_\coop$. Therefore $A\x_B A$ is a double algebra precisely if it satisfies axioms A1, A2, A3 and A8.

\begin{lem}
$\beta_L(A\x_B A)=\Cnt_V(R)\x_B e$ and $\beta_R(A\x_B A)=e\x_B \Cnt_V(L)$.
\end{lem}
\begin{proof}
It suffices to prove the first statement. If $c\in\Cnt_V(R)$ then for $b\in B$
\[
c\ci\Fi_R(b)\oB e=\Fi_R(b)\ci c\oB e=c\star b\oB e=
c\oB b\star e=c\oB e\ci\Fi_L(b).
\]
On the other hand, for all $a_1\x a_2\in A\x_B A$
\begin{gather*}
(a_1\star\Fi_B\Fi_L(a_2))\Fi_R(b)=(a_1\ci\Fi_R(b))\star\Fi_B\Fi_L(a_2)
=a_1\star\Fi_B\Fi_L(a_2\ci\Fi_L(b))\\
=\Fi_R\Fi_B(\Fi_L(a_2)\ci\Fi_L(b))\ci a_1
=\Fi_R\Fi_B\Fi_L(b)\ci\Fi_R\Fi_B\Fi_L(a_2)\ci a_1\\
=\Fi_R(b)\ci(a_1\star\Fi_B(\Fi_L(a_2))
\end{gather*}
where in the second line we used that $\Fi_R\Fi_B$ on $L$ is the antipode inverse, hence antimultiplicative.
\end{proof}

Now we verify the four axioms one-by-one.

\textbf{A1}. It suffices to show that for $c\in \Cnt_V(R)$ and $a_1\x a_2\in A\x_B A$
\[
(c\ci a_1)\x a_2=\beta_B(c\x e)\star (a_1\x a_2)
\]
We compute the RHS,
\begin{gather*}
(c\ci u_k)\star\Fi_B(v_k\star a_1)\x a_2=\Fi_R\Fi_B(v_k\star a_1)\ci c\ci u_k\x a_2\\
=c\ci\Fi_R\Fi_B(v_k\star a_1)\ci u_k\x a_2=c\ci(u_k\star\Fi_B(v_k\star a_1))\x a_2
\end{gather*}
which is the LHS, indeed.

\textbf{A2}. It suffices to show that for $c\in \Cnt_V(R)$ and $a_1\x a_2\in A\x_B A$
\[
(a_1\ci c)\x a_2=\beta_T(c\x e)\star(a_1 \x a_2)
\]
Inserting the definition of $\beta_T$ the RHS reads as
\[
(u_k\ci c)\star\Fi_B(v_k\star a_1)\x a_2=(u_k\star\Fi_B(v_k\star a_1))\ci c\x a_2
=a_1\ci c\x a_2\,.
\]

\textbf{A3}. The left hand side
\[
\beta_B(a_1\x a_2)\star(a'_1\x a'_2)=(a_1\ci u_k)\star\Fi_B((a_2\ci v_k)\star a'_1)\x a'_2
\]
and the right hand side
\[
\beta_L\beta_B(a_1\x a_2)\ci(a'_1\x a'_2)=((a_1\ci u_k)\star\Fi_B\Fi_L(a'_1\ci v_k))\ci a_2\x a'_2
\]
are equal if for all $a_1\x a_2\in A\x_B A$ the map
\[
(a\ci u_k)\star\Fi_B((a_2\ci v_k)\star\under)\colon A\to A
\]
is a right $V$-module map. But this follows from the existence of antipode using
(\ref{eq: SB}) and its the dual basis version in Lemma \ref{lem: S} (2).

\textbf{A8}. Need to show $\beta_T(a_1\x a_2)\star(a'_1\x a'_2)=(a'_1\x a'_2)\ci
\beta_L\beta_T(a_1\x a_2)$.
Inserting the definitions of the $\beta$ we obtain, similarly to the A3 case, that
the map $F=(u_k\ci a_1)\star\Fi_B((v_k\ci a_2)\star\under)\colon A\to A$ need to be left $V$-module map. Using that the antipode is invertible we can compute
\[
F(a')=(u_k\star\Fi_B(v_k\star (a'\ci S^{-1}(a_2))))\ci a_1=a'\ci S^{-1}(a_2)\ci a_1
\]
which is a left $V$-module map, indeed.

\subsection{Hopf algebras} \label{ss: H}

In this subsection we show how Frobenius Hopf algebras can be described as double algebras and point out the difference between the Hopf algebraic and double algebraic antipodes.

\textbf{Notation:} In this subsection $S$ denotes the Hopf algebraic antipode and $\tilde S$ the double algebraic one.

Let $H$ be a Hopf algebra over the commutative ring $k$ and assume the existence of a Frobenius left integral $i\in H$. I.e., $i$ is a left integral and the mapping $f\mapsto f(i)$ is a Frobenius homomorphism on the dual algebra $\hat H$. This is equivalent to that $H$ is a Frobenius $k$-algebra with a Frobenius homomorphism $\lambda\in \hat H$ which is a left integral on $H$. Thus we are in the situation considered in \cite{Pareigis}.
These left integrals are connected by the duality relation $\lambda\la i=1$, or equivalently, $i\la\lambda=\eps$. We will need also the right integrals $\rho=\lambda S^{-1}$ and its dual right integral $j=S(i)$. With $\sigma\in \hat H$ denoting the distinguished grouplike element $\lambda\ra i$ and with $\tau=\sigma S^{-1}$ we have for $a\in H$
\begin{alignat*}{2}
ai&=\eps(a)i&\qquad ja&=j\eps(a)\\
ia&=i\sigma(a)&\qquad aj&=\tau(a)j\\
\lambda\la a&=1\lambda(a)&\qquad a\ra\rho&=1\rho(a)\\
\lambda(i)&=1&\qquad \rho(j)&=1\\
\text{dual basis of $\lambda$: }&i\twoB\o S^{-1}(i\oneB)&\qquad
\text{dual basis of $\rho$: }&i\oneB\o S(i\twoB)
\end{alignat*}
The Nakayama automorphism $\nu$ of $\lambda$ can be computed to be
$\nu=S^{-2}\alpha$ where $\alpha(a)= \sigma\la a$ is an algebra automorphism of $H$. Since the inverse $\alpha^{-1}(a)=a\oneB\sigma S(a\twoB)$, we obtain that $\beta(a):=S\alpha^{-1}S^{-1}(a)=a\ra\sigma$ is another automorphism of $H$. But the coopposite argument shows that also $S\beta^{-1} S^{-1}=\alpha$. This proves $S^2\alpha =\alpha S^2$. Coassociativity implies that $\alpha\beta=\beta\alpha$ and therefore we arrive to the relation
\begin{equation} \label{eq: SabSab}
S^{-1}\alpha\beta S=\beta^{-1}\alpha^{-1}=(\beta\alpha)^{-1}\,.
\end{equation}

A. Connes and H. Moscovici introduced in \cite{Connes-Moscovici} the `deformed' antipode
\begin{equation}
\tilde S:=\nu^{-1}S^{-1}=\alpha^{-1}S=S\beta
\end{equation}
which satisfies
\begin{align}
S(a\star a')&=\tilde S(a')\star S(a)\\
\tilde S(a\star a')&=\tilde S(a')\star \tilde S(a).
\end{align}

\begin{pro}
Let $\bra H,\cdot,1,\Delta,\eps,S\ket$ be a Hopf algebra with a dual pair $(\lambda,i)$ of Frobenius left integrals.
Define the Fourier transforms $\F(a)=a\la \lambda$, $\F^{-1}(\psi)=i\ra\psi S^{-1}$ and the convolution product $a\star a'=\F^{-1}(\F(a)\F(a'))$. Then $\bra H,\cdot,1,\star,i\ket$ is a distributive Frobenius double algebra.
The base algebras are $L=R=k\cdot 1$ and $B=T=k\cdot i$ and
the double algebraic antipode is $\tilde S=\nu^{-1} S^{-1}$.

If $\bra A,\ci,e,\star,i\ket$ is the double algebra arising from a Hopf algebra as above then the Hopf algebra can be reconstructed as $\bra A,\ci,e,\Delta_B,\Fi_B, {\tilde S}^{-1}\nu_L^{-1}\ket$ where $\tilde S$ is the double algebraic antipode and $\nu_L$ is the Nakayama automorphism of $\Fi_L$.
\end{pro}
\begin{proof}
Since $\F$ is a $k$-module isomorphism, $\bra H,\star,i\ket$ is an associative unital algebra. The following alternative expressions for $\star$ will be useful,
\[
a\star a'=\lambda\left(S^{-1}(a'\oneB)\,a\right)\,a'\twoB
=a'\oneB\,\lambda\left(S^{-1}(a'\twoB)\,a\right)
=a\oneB\,\lambda\left(S^{-1}(a')\,a\twoB\right)
\]
At first we compute the base homomorphisms.
\begin{align*}
\Fi_L(a)&=a\star 1=1\lambda(a)\\
\Fi_R(a)&=1\star a=1\rho(a)\\
\Fi_B(a)&=a i=\eps(a)i\\
\Fi_T(a)&=ia=i\sigma(a)
\end{align*}
All of them being scalar multiples of some identity the double algebra axioms
reduce to triviality. Due to the normalization $\eps(1)=\sigma(1)=1$ and $\lambda(i)=1$ and $\rho(i)=\eps(\F^{-1}(\lambda))=\eps(1)=1$ the images of these
$k$-homomorphisms are $k\cdot 1$ and $k\cdot i$, respectively.

The dual bases of $\lambda$ and $\rho$ are also dual bases
of the $\Fi_L$ and $\Fi_R$, respectively. Using the expressions
\[
\eps(a\star a')=\lambda(S^{-1}(a')\,a)\,,\qquad \sigma(a\star a')=\lambda(S^{-1}(\alpha(a'))\,a)
\]
it is easy to check that the $\Fi_B$ and $\Fi_T$ are also Frobenius with dual basis
\[
i\oneB\o i\twoB\,,\qquad \alpha^{-1}(i\oneB)\o i\twoB
\]
respectively. This proves that $\bra A,\cdot,1,\star, i\ket$ is a Frobenius double algebra.

In order to see distributivity we calculate the comultiplications using the dual bases obtained above. Note that sofar $a\oneB\o a\twoB$ stood for the Hopf algebraic coproduct $\Delta(a)$. Fortunately, this is consistent with the double algebraic notation because we find below that $\Delta_B=\Delta$.
\begin{align*}
\Delta_B(a)&=i\oneB\o i\twoB\star a=i\oneB\o \lambda(S^{-1}(a\oneB)i\twoB)\,a\twoB\\
&=a\oneB\o a\twoB\\
\Delta_T(a)&=\alpha^{-1}(i\oneB)\o i\twoB\star a=\alpha^{-1}(i\oneB)\o \lambda(S^{-1}(a\oneB)i\twoB)\,a\twoB\\
&=\alpha^{-1}(a\oneB)\o a\twoB\\
\Delta_L(a)&=ai\twoB\o S^{-1}(i\oneB)\\
\Delta_R(a)&=ai\oneB\o S(i\twoB)
\end{align*}

The distributivity laws:
\begin{align*}
a(a'\star a'')&=aa'\oneB\lambda(S^{-1}(a'')a'\twoB)=
a\oneB a'\oneB\lambda\left(S^{-1}(a'')S^{-1}(a\threeB)a\twoB a'\twoB\right)\\
&=a\oneB a'\o a\twoB a''\\
(a'\star a'')a&=a'\oneB a\lambda(S^{-1}(a'')a'\twoB)=
a'\oneB a\oneB\lambda\left(S^{-1}(a'')a'\twoB a\twoB S(a\threeB)\right)\\
&=a'\oneB a\oneB\lambda\left(S^{-1}(a'' S\nu S(a\threeB))\,a'\twoB a\twoB\right)\\
&=a'a\oneB\star a'' S\nu S(a\twoB)=a' a\oneB\star a''\beta^{-1}(a\twoB)\\
&=a'a\oneT\star a''a\twoT\\
a\star(a'a'')&=a'\oneB a''\oneB \lambda\left(S^{-1}(a''\twoB)S^{-1}(a'\twoB)\,a\right)\\
&=a'\oneB\lambda\left(S^{-1}(a''\twoB)\,ai\twoB\right)\,a''\oneB
\lambda\left(S^{-1}(a''\twoB)S^{-1}(i\oneB)\right)\\
&=(ai\twoB\star a')(S^{-1}(i\oneB)\star a'')
=(a\oneL\star a')(a\twoL\star a'')\\
(a'a'')\star a&=a'\oneB a''\oneB\lambda\left(S^{-1}(a)\,a'\twoB a''\twoB\right)\\
&=a'\oneB\lambda\left(S^{-1}(ai\oneB)\,a'\twoB\right)
a''\oneB\lambda(i\twoB a''\twoB)\\
&=(a'\star a\oneR)(a''\star a\twoR)
\end{align*}

Since antipode exists in distributive double algebras, the next line
\[
\eps(a'\star(a'' a))=\lambda(S^{-1}(a'' a)\, a')=
\lambda(S^{-1}(a'')a' \tilde S(a))=\eps((a'\tilde S(a))\star a'')
\]
proves that the double algebraic antipode is $\tilde S=\nu^{-1}S^{-1}$.

As for the reconstruction of the Hopf algebraic data it suffices to observe that the Nakayama automorphism of $\lambda$ is the Nakayama automorphism of $\Fi_L$.
\end{proof}

The difference between $S$ and $\tilde S$ is a measure of unimodularity. Indeed, if $\sigma=\eps$ then $\nu=S^{-2}$ and therefore $\tilde S=S$.

One may want to check directly that $\cop_T=(\tilde S\o \tilde S)\cop_B^\op$.
As a matter of fact, the familiar identity for the action of the Nakayama automorphism on the dual basis reads for $\lambda$ as
$S^{-1}(i\oneB)\o \nu(i\twoB)=i\twoB\o S^{-1}(i\oneB)$. Therefore
\[
\tilde S(i\twoB)\o\tilde S(i\oneB)=\nu^{-1}S^{-2}(i\oneB)\o i\twoB
=\alpha^{-1}(i\oneB)\o i\twoB
\]
as promised.

\begin{rmk}
Any Hopf algebra $H$ has a Hopf algebroid structure in which the Hopf algebroid antipode is the Hopf algebra antipode \cite{HGD}. In this case $\cop_T=(S\o S)\cop_B^\op$. However, in this Hopf algebroid there exists no 2-sided Frobenius integral, unless $H$ is unimodular.
\end{rmk}

In order to complete the picture we compare the above double algebra of a Hopf algebra to another one which is obtained from the right integral $j$.
Define
\begin{equation}
a\ast a':=\F'^{-1}(\F'(a)\F'(a'))\,\quad a,a'\in H
\end{equation}
where $\F'(a)=\rho\ra a$. Since $\F'=\hat S^{-1}\F S^{-1}$ where $\hat S$ is the antipode of $\hat H$, the new convolution product can be expressed with the old one as
\[
a\ast a'=S(S^{-1}(a')\star S^{-1}(a))\,.
\]
This means that $S$ is an algebra isomorphism $\bra A,\star, i\ket \overset{S}{\longrightarrow} \bra A,\ast,j\ket^\op$.
But then the formula $\tilde S=\nu^{-1}S^{-1}$ and the fact that $\nu$ is a (vertical) algebra isomorphism leads to the isomorphisms
\begin{align}
\bra A,\cdot,1,\star,i\ket&\overset{\nu^{-1}}{\longrightarrow}
\bra A,\cdot,1,\ast,j\ket\\
\bra A,\cdot,1,\star,i\ket&\overset{S}{\longrightarrow}
\bra A,\cdot,1,\ast,j\ket^\op_\coop
\end{align}
of double algebras.

\subsection{Weak Hopf algebras} \label{ss: WHA}

Here we want to generalize the construction of the previous subsection. This weak Hopf algebraic generalization is, however, not straightforward at all because it depends on
the nontrivial theory of invertible modules and half grouplike elements developed by P. Vecserny\'es in \cite{Vecsernyes}. Since that paper
considers finite dimensional weak Hopf algebras over a field $K$, we have to assume the same.

\textbf{Notation:} In this subsection $S$ denotes the weak Hopf algebraic antipode and $\tilde S$ the double algebraic one.

Let $\bra W,\cdot,1,\cop,\eps,S\ket$ be a Frobenius weak Hopf algebra over the field $K$. This means a weak Hopf algebra over $K$ with a left integral $\lambda$ in the dual weak Hopf algebra $\hat W$ which is also a Frobenius homomorphism $W\to K$. Let $i$ denote the dual left integral in $W$, i.e., $\lambda\la i=1$ and $i\la\lambda=\eps$ \cite{BNSz}. We use the notation $W_{L,R}$ for the left/right subalgebras of $W$ defined by the idempotents
\begin{alignat*}{2}
\PL\colon & W\to W_L &\qquad a&\mapsto \eps(1\oneB a)1\twoB\\
\PR\colon & W\to W_R &\qquad a&\mapsto 1\oneB\eps(a 1\twoB)
\end{alignat*}
and $\hat\PL$, $\hat\PR$ stand for the analogue objects for $\hat W$.
The left integral property of $i$ reads as
\begin{equation} \label{eq: WHA FiB}
ai=\PL(a)i\,\qquad a\in W
\end{equation}
but we would like to know what is $ia$? Vecserny\'es proves that $\sigma:=\lambda\ra i$ is a \textit{left grouplike element},
\begin{align*}
\sigma\oneB\o\sigma\twoB&=\eps\oneB\,\sigma\o\eps\twoB\,\sigma
=\hat S(\sigma^{-1})\,\eps\oneB\o\sigma\eps\twoB\\
\hat\PL(\sigma)&=\eps\qquad\hat\PR(\sigma)=\hat S(\sigma)\,\sigma\,.
\end{align*}
Therefore (see p.510 of \cite{Vecsernyes})
\begin{align*}
\lambda\ra ia&=\sigma\ra a=\bra \hat S(\sigma^{-1})\eps\oneB,a\ket\sigma\eps\twoB
=\bra\sigma, \eps\ra (a\ra \hat S(\sigma^{-1}))\ket\\
&=\sigma\ra\PR(a\ra\hat S(\sigma^{-1}))=\lambda\ra i\PR(a\ra\hat S(\sigma^{-1}))
\end{align*}
implying
\begin{equation} \label{eq: WHA FiT}
ia=i\PR(a\ra\hat S(\sigma^{-1})),\qquad a\in W\,.
\end{equation}
Equations (\ref{eq: WHA FiB}) and (\ref{eq: WHA FiT}) are going to define the bottom and top base homomorphisms if $W$ is a double algebra. For that we need a convolution product which we define exactly as in the case of Hopf algebras,
\begin{equation} \label{eq: WHA star}
a\star a':=\F^{-1}(\F(a)\F(a'))\quad\text{where }\F(a):=a\la\lambda\,,
\end{equation}
and can be expressed as
\begin{align}
a\star a'&=a'\ra\hat S^{-1}(a\la\lambda)\notag\\
&=\bra\lambda,S^{-1}(a'\oneB)\,a\ket\, a'\twoB \label{eq: WHA star1}\\
&=a'\fourB S^{-1}(a'\threeB)\,a'\twoB\,\bra\lambda,S^{-1}(a'\oneB)\,a\ket\notag\\
&=a'\threeB\PL(S^{-1}(a'\twoB))\,\bra\lambda,S^{-1}(a'\oneB)\,a\ket\notag\\
&=a'\threeB\PL(S^{-1}(a'\twoB)S^{-1}(\PL(a\oneB))) \,\bra\lambda,S^{-1}(a'\oneB)\,a\twoB\ket\notag\\
&=a'\threeB\PL(S^{-1}(a'\twoB)\,a\oneB) \,\bra\lambda,S^{-1}(a'\oneB)\,a\twoB\ket\notag\\
&=S(S^{-1}(a'\threeB))\,S^{-1}(a'\twoB)\,a\oneB
\,\bra\lambda,S^{-1}(a'\oneB)\,a\twoB\ket\notag\\
&=a\oneB\,\bra\lambda,S^{-1}(a')\,a\twoB\ket \label{eq: WHA star2}
\end{align}
Using these two expressions for the $\star$-product and introducing also the right integral $\rho:=\hat S^{-1}(\lambda)$ we find the base homomorphisms to be
\begin{align}
\Fi_L(a)&=a\star 1=\lambda\la a\\
\Fi_R(a)&=1\star a=a\ra\rho\\
\Fi_B(a)&=ai=\PL(a)i\\
\Fi_T(a)&=ia=i\PR(a\ra\hat S(\sigma^{-1}))
\end{align}
The first two implies that the $L$ and $R$ base ideals are equal to the $W_L$, $W_R$ subalgebras. $B$ is the trivial left $W$-module $Wi$ and $T$ is the space of left integrals in $W$.

Verifying the double algebra axioms the following formula is useful:
\begin{equation} \label{eq: star 4}
(lal')\star(ra'r')=lr(a\star a')l'r'\,,\qquad l,l'\in W_L,\ r,r'\in W_R,\ a,a'\in W
\end{equation}
which follows easily from (\ref{eq: WHA star1}) and (\ref{eq: WHA star2}) using the special form of the comultiplication on $W_L$, $W_R$. It implies that
\[
((\lambda\la a)i)\star a'=(\lambda\la a)(i\star a')=(\lambda\la a)a'
\]
which is \textbf{A1}. Axioms \textbf{A2}, \textbf{A5} and \textbf{A6} can be shown similarly. The remaining four are slightly different. \textbf{A3} and \textbf{A8} follows using that $\lambda\la\under$ is a $W_L$-$W_L$-bimodule map and then (\ref{eq: star 4}). E.g., \textbf{A3} is proven by
\begin{align*}
(\lambda\la(ai))\,a'&=(\lambda\la(\PL(a)i))\,a'=\PL(a)(\lambda\la i)\,a'=\PL(a)\,a'\\
(ai)\star a'&=(\PL(a)i)\star a'=\PL(a)(i\star a')=\PL(a)a'
\end{align*}
The remaining two axioms require the analogous properties of $\rho$.

The $\Fi$'s are all Frobenius homomorphisms as we are going to show now. The dual bases of $\Fi_L$ and $\Fi_R$ are
$i\twoB\o S^{-1}(i\oneB)$ and $i\oneB \o S(i\twoB)$, respectively. The first statement follows from duality of $\lambda$ and $i$, the second from the relation $\Fi_R=S\Fi_L S^{-1}$.
In order to show that the dual basis of $\Fi_B$ is $i\oneB\o i\twoB$ we at first use (\ref{eq: WHA star1}) to calculate
\[
\PL(a\star a')=\bra\lambda,S^{-1}(1\oneB a')\,a\ket\,1\twoB
=S(1\oneB)\,\bra\lambda,S^{-1}(a')\,1\twoB\,a\ket
\]
and then
\begin{align*}
i\oneB\star [\PL(i\twoB\star a)i]&=1\oneB\bra\lambda,S^{-1}(a)1\twoB i\twoB\ket i\oneB
\star i=a\\
[\PL(a\star i\oneB)i]\star i\twoB&=\bra\lambda,S^{-1}(1\oneB i\oneB)\,a\ket 1\twoB i\twoB=a
\end{align*}
as we have claimed. The dual basis of $\Fi_T$ does not follow easily from this because $S$ does not relate them like it did $\Fi_L$ and $\Fi_R$. It is time to introduce
$\tilde S:=\nu^{-1}S^{-1}$ by analogy with the Hopf case, where $\nu$ is the Nakayama automorphism of $\lambda$, therefore \cite{Vecsernyes} $\nu=S^{-2}\alpha$ where $\alpha(a):=\sigma\la a$.
So we have $\tilde S(a)=\sigma^{-1}\la S(a)$, an algebra antiautomorphism of $W$. But it is also an antiautomorphism for the convolution product,
\begin{align*}
\tilde S(a)\star \tilde S(a')&=\bra\lambda,S^{-1}(\tilde S(a')\oneB)\tilde S(a)\ket \tilde S(a')\twoB\\
&=\bra\lambda,S^{-1}(S(a')\oneB)\nu^{-1}S^{-1}(a)\ket\, S(a')\twoB\,\bra\sigma^{-1},S(a')\threeB\ket\\
&=\bra\lambda,a'\threeB\,\nu^{-1}S^{-1}(a)\ket\, S(a'\twoB)\,\bra\sigma^{-1},S(a'\oneB)\ket\\
&=\bra\lambda,S^{-1}(a)\,a'\twoB\ket\,\tilde S(a'\oneB)\\
&=\tilde S(a'\star a)
\end{align*}
and therefore also $\tilde S(i)=i$. Notice that for $r\in W_R$ one has $\tilde S(r)=
\sigma^{-1}\la S(r)=(\sigma^{-1}\la 1)S(r)=S(r)$ because $\sigma^{-1}$ is also left grouplike. Thus we arrive to the desired relation
\begin{align*} 
\Fi_B(\tilde S(a))&=\PL(\sigma^{-1}\la S(a))i=S(\PR(a\twoB))\,i\,\bra\sigma^{-1},
S(a\oneB)\ket=\tilde S(i\,\PR(a\ra\hat S(\sigma^{-1})))\\
&=\tilde S\Fi_T(a)\,.
\end{align*}
Thus $\Fi_T$ is also Frobenius. This finishes the construction of the Frobenius double algebra $\bra W,\cdot,1,\star,i\ket$.

Now we can compute the comultiplications.
\begin{align} \label{eq: WHA DeltaB}
\cop_B(a)&=i\oneB\oB i\twoB\star a=i\oneB\oB \bra\lambda,S^{-1}(a\oneB)i\twoB\ket a\twoB \notag\\
&=a\oneB\oB a\twoB
\end{align}
so we can keep the notation $a\oneB\o a\twoB$ for $\cop_B$ but putting the $B$-module tensor product instead. The precise relation is, of course, that $\cop_B=\Pi\cop$ where $\Pi\colon A\o A\twoheadrightarrow A\oB A$ is the canonical epimorphism. About
\begin{equation} \label{eq: WHA DeltaL}
\cop_L(a)=ai\twoB\oL S^{-1}(i\oneB)
\end{equation}
is nothing to say. The remaining comultiplications will not be used explicitly, so it suffices to remark that because of $\tilde S\Fi_L=\Fi_R\tilde S$ the right dual basis and therefore the right comultiplication $\cop_R$, too, can be calculated from those of $\Fi_L$.

What we already know about $\tilde S$ is at the half-way of proving that $\tilde S$ is the antipode of this double algebra. We have seen that $\tilde S$ is a double algebra isomorphism $W\to W^\op_\coop$. So it suffices to verify the antipode axioms (\ref{eq: SB}) and (\ref{eq: SL}).
\begin{align*}
\Fi_B((a'\tilde S(a))\star a'')&=S(1\oneB)\,\bra\lambda,S^{-1}(a'')1\twoB a'\,\nu^{-1}S^{-1}(a)\ket\,i\\
&=S(1\oneB)\,\bra\lambda,S^{-1}(a''a)1\twoB a'\ket\,i\\
&=\Fi_B(a'\star(a''a))\\
\Fi_L(a'(a''\star\tilde S(a)))&=\bra\lambda,S^{-1}(a\oneB)\,a'(a''\star\tilde S(a))\ket 1\twoB\\
&=\bra\lambda,S^{-1}(a\oneB)\,a'\bra\lambda,S^{-1}(S(a)\oneB)\,a''\ket
(\sigma^{-1}\la S(a)\twoB)\ket\,1\twoB\\
&=\bra\lambda,a\twoB a''\ket\, \bra\lambda,S^{-1}(1\oneB)\,a'\nu^{-1}S^{-1}(a\oneB)\ket\,1\twoB\\
&=\bra\lambda,a\twoB a''\ket\,\bra\lambda,S^{-1}(1\oneB a\oneB)\,a'\ket\,1\twoB\\
&=\bra\lambda,S^{-1}(1\oneB)\,a\twoB a''\ket\,\bra\lambda,S^{-1}(a\oneB)\,a'\ket\,1\twoB\\
&=\bra\lambda,S^{-1}(1\oneB)(a'\star a)a''\ket\,1\twoB\\
&=\Fi_L((a'\star a)a'')
\end{align*}

Having the antipode distributivity is equivalent to that two comultiplications, $\cop_B$ and $\cop_L$ for example, are multiplicative. Multiplicativity of the others will result from applying $\tilde S$. Multiplicativity of $\cop_B$ is an obvious consequence of multiplicativity of the weak Hopf algebraic comultiplication because of (\ref{eq: WHA DeltaB}). Multiplicativity of $\cop_L$ is proven by the calculation
\begin{align*}
\cop_L(a)\star\cop_L(a')&=ai\twoB\star a'i_{(2')}\oL S^{-1}(i\oneB)\star S^{-1}(i_{(1')})\\
&=\bra\lambda,S^{-1}(a'\oneB i_{(2')})\,a i\twoB\ket\,a'\twoB i_{(3')}\,\oL\,
S^{-1}(i\oneB)\star S^{-1}(i_{(1')})\\
&=a'\twoB i_{(3')}\,\oL\, S^{-1}(a'\oneB i_{(2')})\,a\star S^{-1}(i_{(1')})\\
&=a'\twoB i\fourB\,\oL\,\bra\lambda, S^{-2}(i\twoB)S^{-1}(i\threeB)S^{-1}(a'\oneB)\,a\ket\,S^{-1}(i\oneB)\\
&=a'\twoB\PL(i\twoB) i\threeB\,\oL\,\lambda, S^{-1}(a'\oneB)\, a\ket\,S^{-1}(i\oneB)\\
&=(a\star a')i\twoB\,\oL\, S^{-1}(i\oneB)\ = \ \cop_L(a\star a')\,.
\end{align*}
Let us summarize what we have proven above.
\begin{pro}
Let $\bra W,\cdot,1,\cop,\eps\ket $ be a finite dimensional weak Hopf algebra over the field $K$ and let $i\in W$ be a nondegenerate left integral. Then $\bra W,\cdot,1,\star,i\ket$ is a distributive Frobenius double algebra where the convolution $\star$ is defined by (\ref{eq: WHA star}).
\end{pro}

Reconstruction of the weak Hopf algebra from its double algebra is not possible completely. Of course, the antipode $S$ is completely reconstructed from $\tilde S$ and from the Nakayama automorphism of $\Fi_L$. From the weak bialgebra structure, however, only two bialgebroids remain, $V_B$ and $V_T$, in the double algebra. The restriction $\eps|_L$ of the counit cannot be reconstructed. This means that the double algebra has to be supplied with the data of an index one Frobenius functional $\psi\colon L\to K$ on the separable $K$-algebra $L$ \cite[Propositions 7.3, 7.4]{KSz}.

\subsection{Hopf algebroids} \label{ss: HGD}

Here we will show that a Hopf algebroid in the sense of \cite{HGD} with a chosen two-sided Frobenius integral $i$ is a distributive Frobenius double algebra, thereby finishing the proof of Theorem \ref{thm: main}. This task is much simpler than the analogous ones in case of Hopf and weak Hopf algebras because the
Hopf algebroid antipode is flexible enough to become a double algebraic antipode. Therefore we do not need distinguished grouplike elements to deform the antipode with.

Let $\A=\bra A,B,s_L,t_L,\gamma_L,\pi_L\ket$ be a left bialgebroid over $B$ in the category of $k$-modules and let $\bra \A,S\ket$ be a Hopf algebroid with a Frobenius (called `nondegenerate' in \cite{HGD}) left integral $i\in A$. Without loss of generality we may
assume that $i=S(i)$ by deforming the antipode if necessary \cite[Proposition 5.13]{HGD}. In the sequel we use the double algebraic notation for tensor products like
$\oB$, $\oL$,\dots cf. the dictionary at the end.
Let $\phi\colon \,_BA_B\to \,_BB_B$ be the unique Frobenius homomorphisms
with dual basis $i\twoT\oL S^{-1}(i\oneT)$. Let $_*A=\Hom(_BA,\,_BB)$ and define
the Fourier transformations
\begin{alignat*}{2}
\F\colon&A\to\,_*A&\qquad a&\mapsto\phi(\under a)\\
\F^{-1}\colon&_*A\to A&\qquad f&\mapsto i\twoT\, s_L f(S^{-1}(i\oneT))
\end{alignat*}
and the convolution product
\[
a\star a':=\F^{-1}(\F(a)\bullet\F(a'))
\]
where the dot is the multiplication on $_*A$ that is the opposite of the one given in \cite[Eq. (42)]{KSz}. (This oppositeness is to comply with equations (\ref{eq: clue1}), (\ref{eq: clue2}).)
Now it is obvious that $\star$ is associative with unit $i$.
A concrete formula can be calculated as
\begin{align}
a\star a'&=i\twoT\,s_L
\phi\left([S^{-1}(i\oneT)\oneB\cdot\phi(S^{-1}(i\oneT)\twoB\,a')]\,a\right)\notag\\
&=i\threeT\,s_L\phi\left(S^{-1}s_L\phi(S^{-1}(i\oneT)\,a')\,S^{-1}(i\twoT)\,a\right)
\notag\\
&={a'}\twoT i\threeT\,
s_L\phi\left(S^{-1}({a'}\oneT i\twoT\,s_L\phi(S^{-1}(i\oneT)))\,a\right) \notag\\
&={a'}\twoT\,s_L\phi(S^{-1}({a'}\oneT)\,a) \label{eq: HGD star1}
\end{align}

So we have the base homomorphisms
\begin{align}
\Fi_L(a)&=a\star 1=s_L\phi(a)\\
\Fi_R(a)&=1\star a=a\twoT\,s_L\phi(S^{-1}(a\oneT))\\
\Fi_B(a)&=ai=s_L\PL(a)\,i\\
\Fi_T(a)&=ia=i\,s_R\PR(a)
\end{align}
From the above formula for $a\star a'$ it is immediate that
\[
(lal')\star (ra'r')=lr(a\star a')l'r'\,\quad l,l'\in s_L(L),\ r,r'\in s_R(R).
\]
Therefore axioms \textbf{A1}, \textbf{A2}, \textbf{A5} and \textbf{A6} can be shown easily like in the weak Hopf case. The proof of the remaining ones goes as follows.
\begin{align*}
\Fi_L\Fi_B(a)a'&=s_L\phi(s_L\PL(a)i)a'=s_L\PL(a)a'=(s_L\PL(a)i)\star a'\\
&=\Fi_B(a)\star a'\\
\Fi_R\Fi_B(a)a'&=i\twoT\,s_L\phi(S^{-1}(i\oneT)S^{-1}(s_L\PL(a)))a'
=S^{-1}(s_L\PL(a))a'\\
&=a'\star (S^{-1}(s_L\PL(a))i)=a'\star (s_L\PL(a)i)=a'\star\Fi_B(a)\\
\end{align*}
\begin{align*}
a\Fi_R\Fi_T(a')&=ai\twoT\,s_R\PR(a')\,s_L\phi(S^{-1}(i\oneT))=as_R\PR(a')a\star (is_R\PR(a'))\\
&=a\star\Fi_T(a')\\
a\Fi_L\Fi_T(a')&=a\,s_L\phi(is_R\PR(a'))=as_L\phi(iS^{-1}(s_R\PR(a')))\\
&=as_L\PL S^{-1}(s_R\PR(a'))=aS^{-1}s_R\PR(a') =(iS^{-1}s_R\PR(a'))\star a\\
&=(is_R\PR(a'))\star a=\Fi_T(a')\star a
\end{align*}
Thus we have a double algebra.

Next we want to show that $S$ is a $\star$-antiautomorphism.
For that we need an alternative formula for the convolution. At first, notice that $S(i)=i$ implies
\begin{equation} \label{eq: L duba}
i\twoT\oL S^{-1}(i\oneT)=S(i\oneB)\oL i\twoB
\end{equation}
Using also the calculation
\begin{align}
a\star i\oneB \oB i\twoB&={i\oneB}\twoT\,
s_L\phi(S^{-1}({i\oneB}\oneT)\,a)\oB i\twoB\notag\\
&={i\twoT}\oneB\,s_L\phi(S^{-1}(i\oneT)\,a)\oB {i\twoT}\twoB\notag\\
&=\gamma_L\left(i\twoT\,s_L\phi(S^{-1}(i\oneT)\,a)\right)\notag\\
&=a\oneB \oB a\twoB \label{eq: a*ioi}
\end{align}
we obtain
\begin{align}
a\star a'&=a\star S^{-1}\left(S(i\oneB)\,s_L\phi(i\twoB S(a'))\right) \notag\\
&=a\star t_L\phi(i\twoB S(a'))\,i\oneB
=t_L\phi(i\twoB S(a'))\,(a\star i\oneB)\notag\\
&=t_L\phi(a\twoB S(a'))\, a\oneB \label{eq: HGD star2}
\end{align}
A consequence is that $S$ is $\star$-antimultiplicative,
\begin{align}
S^{-1}(a\star a')&=t_L\phi(S^{-1}({a'}\oneT)\,a)\,S^{-1}({a'}\twoT)
=t_L\phi(S^{-1}(a')\twoB\,S(S^{-1}(a)))\,S^{-1}(a')\oneB\notag\\
&=S^{-1}(a')\star S^{-1}(a)
\end{align}

\begin{table}
\begin{tabular}{| >{$}c<{$}| >{$}c<{$}|}
\hline
\text{Hopf algebroid}&\text{double algebra}\\ \hline
1 & e \\
L & B\\
R & T \\
s_L & \Fi_L|_B\\
t_L & \Fi_R|_B\\
s_R & \Fi_R|_T\\
t_R & \Fi_L|_T\\
s_L(L)=t_R(R) & L\\
t_L(L)=s_R(R) & R\\
\pi_L & \Fi_B\\
\pi_R & \Fi_T\\
_LA_L & _BA_B\\
^RA^R & _TA_T\\
_LA^L & _LA_L\\
_RA^R & _RA_R\\
\gamma_L(a)=a\oneB\oL a\twoB & \Delta_B(a)=a\oneB\oB a\twoB\\
\gamma_R(a)=a\oneT\oR a\twoT & \Delta_T(a)=a\oneT\oT a\twoT\\
\text{a two sided Frobenius integral }i& i\\
S(i\oneB)\oL i\twoB a&\Delta_L(a)=a\oneL\oL a\twoL\\
ai\oneT\oR S(i\twoT)&\Delta_R(a)=a\oneR\oR a\twoR\\
S & S\\ \hline
\end{tabular}
\vskip 0.5truecm
\caption{The dictionary}
\label{tab: dictionary}
\end{table}

By construction the $\Fi_L$ is a Frobenius homomorphism to $L:=s_L(B)\subset A$ with dual basis (\ref{eq: L duba}). From double algebraic experience one conjectures that $\Fi_B$ has dual basis $i\oneB\oB i\twoB$. Indeed,
\begin{gather*}
s_L\PL(a\star i\oneB)i\star i\twoB=s_L\PL(t_L\phi(a\twoB S(i\oneB))\,a\oneB)\,i\twoB\\
=s_L\PL(a\oneB)\,s_L\phi(a\twoB S(i\oneB))\,i\twoB=s_L\PL(a\oneB)\,a\twoB=a\\
i\oneB\star (s_L\PL(i\twoB\star a)i)=i\oneB\star S^{-1}s_L\PL(i\twoB\star a)i\\
=t_L\PL(t_L\phi(i\threeB S(a))\,i\twoB)\,i\oneB
=t_L\left(\PL(i\twoB)\,\phi(i\threeB S(a))\right)\,i\oneB\\
=t_L\phi(i\twoB S(a))\,i\oneB=a
\end{gather*}
Using also $\Fi_R=S^{-1}\Fi_L S$ and $\Fi_T=S^{-1}\Fi_B S$ they also are Frobenius homomorphisms. Summarizing, we found the dual bases
\begin{align*}
x_j\oL y_j&=S(i\oneB)\oL i\twoB\\
u_k\oB v_k&=i\oneB\oB i\twoB\\
x^j\oR y^j&=i\oneT\oR S(i\twoT)\\
u^k\oT v^k&=i\oneT\oT i\twoT
\end{align*}
Together with equation (\ref{eq: a*ioi}) the 2$^{\text{nd}}$ and 4$^{\text{th}}$ of these imply $\cop_B=\gamma_L$ and $\cop_T=\gamma_R$. In particular $\cop_B$ is multiplicative. But we want to prove distributivity directly. Since $S$ is a double algebra antiautomorphism, it suffices to verify (\ref{eq: DB}) and (\ref{eq: DR}).
\begin{align*}
(a\oneB a')\star (a\twoB a'')&={a\twoB}\twoT {a''}\twoT\,
s_L\phi\left(S^{-1}({a\twoB}\oneT {a''}\oneT)\,a\oneB a'\right)\\
&=a\twoT{a''}\twoT\,s_L\phi\left(S^{-1}({a''}\oneT)S^{-1}({a\oneT}\twoB){a\oneT}\oneB a'
\right)\\
&=a\twoT t_R\PR(a\oneT){a''}\twoT\,s_L\phi(S^{-1}({a''}\oneT)\,a')\\
&=a(a'\star a'')
\end{align*}
where in the 3$^{\text{rd}}$ row we used the Takeuchi property for $\cop_T(a'')$.
\begin{align*}
(a'\star a\oneR)(a''\star a\twoR)&=(a'\star(ai\oneT))(a''\star S(i\twoT))\\
&=(a'\star(ai\oneT))\,t_L\phi(a''\twoB S^2(i\twoT))\,a''\oneB\\
&=\left(a'\star a\,S^{-1}\left(s_L\phi(a''\twoB i\twoT)\,S^{-1}(i\oneT)\right)\right)
\,a''\oneB\\
&=(a'\star aS^{-1}(a''\twoB))\,a''\oneB\\
&=t_L\phi(a'\twoB a''\twoB S(a))\,a'\oneB a''\oneB\ =\ (a'a'')\star a
\end{align*}
Thus we have proven that $\bra A,\cdot,1,\star,i\ket$ is a distributive Frobenius double algebra. For then the antipode exists $S$ will be proven to be the antipode
once we verify one of its defining relations, let us say the last one in Lemma \ref{lem: S} (2) which says
\[
S(a)i\oneT\oT i\twoT  = i\oneT\oT a\,i\twoT\,.
\]
But this is precisely the left integral property of $i$, cf. \cite[Lemma 5.2]{HGD}.
The dictionary in Table \ref{tab: dictionary} helps to compare the Hopf algebroid notations of \cite{HGD} and \cite{FrobD2} with the double algebraic notations.

This finishes the proof of Theorem \ref{thm: main}. The calculations presented in this subsection perhaps illustrate the advantage of the double algebraic view as opposed to the bialgebroid view of Hopf algebroids, at least when a Frobenius integral is present.
The question is still pending whether anything remains from the double algebra structure in the absence of good integrals?

\end{document}